\documentclass[12pt]{article}
\usepackage{amssymb}

\title{On the concepts of intertwining operator and tensor
product module in vertex operator algebra theory}
\author{Yi-Zhi Huang, James Lepowsky, Haisheng Li and Lin Zhang}
\date{}

\def \Z{\mathbb Z}
\def \C{\mathbb C}

\def \N{\mathbb N}
\def \H{{\mathcal{H}}}

\def \Y{{\mathcal{Y}}}

\def \swt{{\rm \scriptstyle wt \,}}
\def \wt{{\rm wt \,}}

\def \Res{{\rm Res \,}}
\def \End{{\rm End \,}}

\def \Hom{{\rm Hom \,}}

\def \<{\langle}
\def \>{\rangle}

\def \g{\mathfrak{g}}

\def \bconj{\begin{conj}\label}
\def \econj{\end{conj}}
\def \be{\begin{equation}\label}
\def \ee{\end{equation}}
\def \bex{\begin{exa}\label}
\def \eex{\end{exa}}

\def \bl{\begin{lem}\label}
\def \el{\end{lem}}
\def \bt{\begin{thm}\label}
\def \et{\end{thm}}
\def \bp{\begin{prop}\label}
\def \ep{\end{prop}}
\def \br{\begin{rem}\label}
\def \er{\end{rem}}
\def \bc{\begin{coro}\label}
\def \ec{\end{coro}}
\def \bd{\begin{de}\label}
\def \ed{\end{de}}
\def \pf{{\bf Proof. }}

\newtheorem{thm}{Theorem}[section]
\newtheorem{prop}[thm]{Proposition}
\newtheorem{coro}[thm]{Corollary}
\newtheorem{conj}[thm]{Conjecture}
\newtheorem{exa}[thm]{Example}
\newtheorem{lem}[thm]{Lemma}
\newtheorem{rem}[thm]{Remark}
\newtheorem{de}[thm]{Definition}

\def \nno{\nonumber}

\def \pf{{\it Proof}\hspace{2ex}}
\def \epf{\hspace{1em}$\square$}
\def \epfv{\hspace{1em}$\square$\vspace{1em}}

\newtheorem{propo}[thm]{Proposition}
\newtheorem{corol}[thm]{Corollary}
\newtheorem{exam}[thm]{Example}

\newtheorem{rema}[thm]{Remark}
\newtheorem{defi}[thm]{Definition}

 \makeatletter
\newlength{\@pxlwd} \newlength{\@rulewd} \newlength{\@pxlht}
\catcode`.=\active \catcode`B=\active \catcode`:=\active \catcode`|=\active
\def\sprite#1(#2,#3)[#4,#5]{
   \edef\@sprbox{\expandafter\@cdr\string#1\@nil @box}
   \expandafter\newsavebox\csname\@sprbox\endcsname
   \edef#1{\expandafter\usebox\csname\@sprbox\endcsname}
   \expandafter\setbox\csname\@sprbox\endcsname =\hbox\bgroup
   \vbox\bgroup
  \catcode`.=\active\catcode`B=\active\catcode`:=\active\catcode`|=\active
      \@pxlwd=#4 \divide\@pxlwd by #3 \@rulewd=\@pxlwd
      \@pxlht=#5 \divide\@pxlht by #2
      \def .{\hskip \@pxlwd \ignorespaces}
      \def B{\@ifnextchar B{\advance\@rulewd by \@pxlwd}{\vrule
         height \@pxlht width \@rulewd depth 0 pt \@rulewd=\@pxlwd}}
      \def :{\hbox\bgroup\vrule height \@pxlht width 0pt depth
0pt\ignorespaces}
      \def |{\vrule height \@pxlht width 0pt depth 0pt\egroup
         \prevdepth= -1000 pt}
   }
\def\endsprite{\egroup\egroup}
\catcode`.=12 \catcode`B=11 \catcode`:=12 \catcode`|=12\relax
\makeatother

\def\hboxtr{\FormOfHboxtr} 
\sprite{\FormOfHboxtr}(25,25)[0.5 em, 1.2 ex] 
:BBBBBBBBBBBBBBBBBBBBBBBBB |
:BB......................B |
:B.B.....................B |
:B..B....................B |
:B...B...................B |
:B....B..................B |
:B.....B.................B |
:B......B................B |
:B.......B...............B |
:B........B..............B |
:B.........B.............B |
:B..........B............B |
:B...........B...........B |
:B............B..........B |
:B.............B.........B |
:B..............B........B |
:B...............B.......B |
:B................B......B |
:B.................B.....B |
:B..................B....B |
:B...................B...B |
:B....................B..B |
:B.....................B.B |
:B......................BB |
:BBBBBBBBBBBBBBBBBBBBBBBBB |
\endsprite

\sprite{\FormOfShboxtr}(25,25)[0.3 em, 0.72 ex] 
:BBBBBBBBBBBBBBBBBBBBBBBBB |
:BB......................B |
:B.B.....................B |
:B..B....................B |
:B...B...................B |
:B....B..................B |
:B.....B.................B |
:B......B................B |
:B.......B...............B |
:B........B..............B |
:B.........B.............B |
:B..........B............B |
:B...........B...........B |
:B............B..........B |
:B.............B.........B |
:B..............B........B |
:B...............B.......B |
:B................B......B |
:B.................B.....B |
:B..................B....B |
:B...................B...B |
:B....................B..B |
:B.....................B.B |
:B......................BB |
:BBBBBBBBBBBBBBBBBBBBBBBBB |
\endsprite

\newsavebox{\BoxTimes}\sbox{\BoxTimes}{$\boxtimes$}
\def\warning{\kern 2pt \raisebox{-0.5 pt}{\FormOfWarning}\kern 3 pt}
\sprite{\FormOfWarning}(25,25)[0.75em, 0.75em]
:BBBBBBBBBBBBBBBBBBBBBBBBB |
:B.......................B |
:B.......................B |
:B.......................B |
:B.......................B |
:B....BBBBBBBBBBBBBB.....B |
:B....BBBBBBBBBBBBBBBB...B |
:B....BBBBBBBBBBBBBBBB...B |
:B...............BBBBB...B |
:B..............BBBBB....B |
:B.............BBBBB.....B |
:B...........BBBBB.......B |
:B.........BBBBB.........B |
:B.......BBBBB...........B |
:B.....BBBBB.............B |
:B....BBBBB..............B |
:B....BBBB...............B |
:B....BBBBBBBBBBBBBBB....B |
:B....BBBBBBBBBBBBBBB....B |
:B......BBBBBBBBBBBBB....B |
:B.......................B |
:B.......................B |
:B.......................B |
:B.......................B |
:BBBBBBBBBBBBBBBBBBBBBBBBB |
\endsprite

\begin{document}
\maketitle

\begin{abstract}
We produce counterexamples to show that in the definition of the
notion of intertwining operator for modules for a vertex operator
algebra, the commutator formula cannot in general be used as a
replacement axiom for the Jacobi identity.  We further give a
sufficient condition for the commutator formula to imply the Jacobi
identity in this definition.  Using these results we illuminate the
crucial role of the condition called the ``compatibility condition''
in the construction of the tensor product module in vertex operator
algebra theory, as carried out in work of Huang and Lepowsky.  In
particular, we prove by means of suitable counterexamples that the
compatibility condition was indeed needed in this theory.
\end{abstract}

\renewcommand{\theequation}{\thesection.\arabic{equation}}
\renewcommand{\thethm}{\thesection.\arabic{thm}}
\setcounter{equation}{0}
\setcounter{thm}{0}

\section{Introduction}

With motivation from both mathematics and physics, a tensor product
theory of modules for a vertex operator algebra was developed in
\cite{tensor1}--\cite{tensor3} and \cite{tensor4}. This theory has had
a number of applications and has been generalized to additional
important settings. The main purpose of this paper is to study and
elucidate certain subtle aspects of this theory and of the fundamental
notion of intertwining operator in vertex operator algebra theory. In
particular, we answer a number of questions that have arisen in the
theory.

A central theme in vertex (operator) algebra theory is that the theory
cannot be reduced to Lie algebra theory, even though the theory is,
and always has been, intimately related to Lie algebra theory.  Before
we discuss the notions of intertwining operator and tensor product
module in vertex operator algebra theory, we first recall that for the
three basic notions of vertex (operator) algebra, of module and of
intertwining operator, there is a uniform main axiom: the Jacobi
identity; see \cite{FLM} and \cite{FHL}.  For the notion of vertex
(operator) algebra itself, the ``commutator formula'' or
``commutativity'' or ``weak commutativity'' (``locality'') can
alternatively be taken to be the main axiom; see \cite{FLM},
\cite{FHL}, \cite{DL}, \cite{Li-local}; cf.\ \cite{LL}.  (Borcherds's
original definition of the notion of vertex algebra \cite{B} used
``skew-symmetry'' and the ``associator formula''; cf.\ \cite{LL}.)  By
contrast, for the notion of module for a vertex (operator) algebra,
the ``associator formula'' or ``associativity'' or ``weak
associativity'' can be used in place of the Jacobi identity as the
main axiom; see \cite{B}, \cite{FLM}, \cite{FHL}, \cite{DL},
\cite{Li-local}; cf.\ \cite{LL}.  And furthermore, in the definition of
this notion of module, the commutator formula (or commutativity or
weak commutativity or locality) {\it cannot} be taken as the main
replacement axiom.  This is actually an easy observation, and in an
Appendix of the present paper we give examples to verify and
illustrate this fact.  This already illustrates how vertex operator
algebra theory cannot be reduced to Lie algebra theory.

In the main text of the present paper, we discuss and analyze the
extent to which the Jacobi identity can be replaced by the ``commutator
formula'' in the definition of the notion of intertwining operator among
modules for a vertex (operator) algebra. This will in particular
explain the crucial nature of the ``compatibility condition'' in
\cite{tensor1}--\cite{tensor3}. (We shall recall in Section
\ref{secquasi} below the precise meaning of the term ``commutator
formula'' in the context of intertwining operators.)

In this paper we assume the reader is familiar with the basic concepts
in the theory of vertex operator algebras, including modules and
intertwining operators; we shall use the theory as developed in
\cite{FLM}, \cite{FHL} and \cite{LL}, and the terminology and notation
of these works.

In the tensor product theory of modules for a vertex operator algebra
(see \cite{tensor1}--\cite{tensor3}), the tensor product functor
depends on an element of a certain moduli space of three-punctured
spheres with local coordinates at the punctures.  In this paper we
shall focus on the important moduli space element denoted $P(z)$ in
\cite{Hbook} and \cite{tensor1}, where $z$ is a nonzero complex
number.  The corresponding tensor product functor is denoted
$W_1\boxtimes_{P(z)} W_2$ for modules $W_1$ and $W_2$ for a suitable
vertex operator algebra $V$. This tensor product module
$W_1\boxtimes_{P(z)} W_2$ can be constructed by means of its
contragredient module, which in turn can be realized as a certain
subspace $W_1\hboxtr_{P(z)}W_2$ of $(W_1\otimes W_2)^*$ (the dual
vector space of the vector space tensor product $W_1\otimes W_2$ of
$W_1$ and $W_2$). The elements of $W_1\hboxtr_{P(z)}W_2$ satisfy the
``lower truncation condition'' and the ``$P(z)$-compatibility
condition'' defined and discussed in
\cite{tensor1}--\cite{tensor3}. It was proved in
\cite{tensor1}--\cite{tensor3} that these two conditions together
imply the Jacobi identity, and hence that any element of $(W_1\otimes
W_2)^*$ satisfying these two conditions generates a weak module for
the vertex operator algebra $V$. We show here that the converse of
this statement is not true in general; specifically, an element of
$(W_1\otimes W_2)^*$ generating a weak $V$-module does not need to
satisfy the compatibility condition. It follows in particular that the
largest weak $V$-module in $(W_1\otimes W_2)^*$, which we shall write
as $W_1\warning W_2$ and shall read as ``$W_1$ warning $W_2$,'' can
indeed be (strictly) larger than the desired space,
$W_1\hboxtr_{P(z)}W_2$. In particular, for each of the examples, or
really counterexamples, that we give, we will see that when the
modules $W_1$ and $W_2$ are taken to be $V$ itself, neglecting the
compatibility condition results in a $V$-module $V\warning V$ whose
contragredient module is indeed (strictly) larger (in the sense of
homogeneous subspaces) than the correct tensor product
$V\boxtimes_{P(z)} V$, which is naturally isomorphic to $V$
itself. This of course shows that the compatibility condition cannot
in general be omitted.

In Section 2 of this paper we will show that the compatibility
condition for elements of $(W_1\otimes W_2)^*$ reflects in a precise
way the Jacobi identity for intertwining operators and intertwining
maps, while on the other hand, the {\it Jacobi identity} (which, as we
have been discussing, is implied by the compatibility condition) for
elements of $(W_1\otimes W_2)^*$ reflects in a precise way the {\it
commutator formula} for intertwining operators and intertwining maps.

Thus we have the natural question (which we already mentioned above):
In the notion of intertwining operator, can the commutator formula be
used as a replacement axiom for the Jacobi identity?  In other words,
does the Jacobi identity imply the compatibility condition in
$(W_1\otimes W_2)^*$?  As one should expect, the answer is no.  We
shall correspondingly call a ``quasi-intertwining operator'' an
operator satisfying all the conditions for an intertwining operator
except that the Jacobi identity in the definition is replaced by the
commutator formula. We shall exhibit a straightforward counterexample
(a quasi-intertwining operator that is not an intertwining operator)
when the vertex algebra is constructed from a commutative associative
algebra with identity.  However, when the vertex algebra has a
conformal vector and nonzero central charge (for instance, when the
vertex algebra is a vertex operator algebra with nonzero central
charge), we will see that the answer is instead yes---the commutator
formula indeed implies the Jacobi identity in this case.  We establish
this and related results and construct relevant counterexamples in
Sections 3 and 4.

As we also show, all these results actually hold in the presence of
logarithmic variables, when the modules involved are only direct sums
of generalized $L(0)$-eigenspaces, instead of $L(0)$-eigenspaces (see
\cite{Mi}, \cite{HLZ1}, \cite{HLZ2} for these notions).

In the Appendix we show that unless a vertex (operator) algebra is
one-dimensional, there exists a non-module that satisfies all the
axioms for a module except that the Jacobi identity is replaced by the
commutator formula.

We would like to add a few more words concerning why one should expect
that for the notion of intertwining operator, the commutator formula
does not imply the Jacobi identity (and consequently, the Jacobi
identity does not imply the compatibility condition).  Consider the
elementary situation in which a vertex algebra $V$ is based on a
commutative associative algebra.  A $V$-module is exactly the same as
a module for the underlying commutative associative algebra, since the
notion of $V$-module can be described via associativity.  In
particular, if $\dim V > 1$, a $V$-module is not in general the same
as a module for $V$ viewed as a commutative Lie algebra (since this
would amount to a vector space of commuting operators acting on the
module); commutativity cannot be used as a replacement axiom in the
definition of the notion of module.  Thus in this situation, the
notion of $V$-module is essentially ring-theoretic and not
Lie-algebra-theoretic, while the notion of quasi-intertwining operator,
based as it is on the commutator formula, is essentially
Lie-algebra-theoretic and not ring-theoretic.  These considerations
motivated our (straightforward) construction of examples showing that
the commutator formula does not imply the Jacobi identity in the
definition of the notion of intertwining operator.

\paragraph{Acknowledgments}
Y.-Z.H., J.L. and L.Z. gratefully acknowledge partial support from
NSF grant DMS-0070800 and H.L. gratefully acknowledges partial support from
an NSA grant.

\setcounter{equation}{0}
\setcounter{thm}{0}
\section{Quasi-intertwining operators and the compatibility condition}
\label{secquasi}

Throughout this section we let $(V,Y,{\bf 1},\omega)$ be a vertex
operator algebra (in the precise sense of \cite{FHL}, \cite{FLM},
\cite{LL} or \cite{tensor1}--\cite{tensor3}).  (Recall that
$V=\coprod_{n\in \Z}V_{(n)}$ is the underlying $\Z$-graded vector space,
$Y$ is the vertex operator map, ${\bf 1}$ is the
vacuum vector, $\omega$ is the conformal vector, and
$Y(\omega,x)=\sum_{n\in\Z}L(n)x^{-n-2}$.)  Let $z$ be a
fixed nonzero complex number.  In this section we first define the
notion of quasi-intertwining operator and
quasi-($P(z)\mbox{-}$)intertwining map, generalizing (and weakening)
the notions of intertwining operator (see \cite{FHL}) and of
($P(z)$-)intertwining map (see Section 4 of \cite{tensor1}). We
establish the correspondence between these two notions, similar to the
correspondence between intertwining operators and intertwining
maps. We then show that for $V$-modules $W_1$ and $W_2$, a
quasi-$P(z)$-intertwining map of type ${W_3\choose W_1\,W_2}$ gives a
weak $V$-module inside $(W_1\otimes W_2)^*$. The notion of logarithmic
quasi-intertwining operator is also defined in this section.

The notion of quasi-intertwining operator is defined in the same way
as the notion of intertwining operator except that the Jacobi identity
is replaced by the ``commutator formula.''  We now in fact give the
definition in the greater generality of weak $V$-modules; a {\it weak
module} for our vertex operator algebra $V$ is a module for $V$ viewed
as a vertex algebra, in the sense of Definition 4.1.1 in \cite{LL}.
\begin{defi}{\rm
Let $(W_1, Y_{1})$, $(W_2, Y_{2})$ and $(W_3, Y_{3})$
be weak $V$-modules. A {\em
quasi-inter\-twining operator of type ${W_3\choose W_1\, W_2}$} is a
linear map ${\cal Y}: W_1\otimes W_2\to W_3\{x\}$ (the space of formal
series in complex powers of $x$ with coefficients in $W_3$), or
equivalently,
\begin{eqnarray*}
W_1 &\to & (\Hom(W_2,W_3))\{x\}\nno\\
w & \mapsto & {\cal Y}(w,x) =\sum_{n\in\C}w_nx^{-n-1}\;\;\;
(\mbox{where}\;\; w_n \in \Hom(W_2,W_3))
\end{eqnarray*}
such that for $v \in V$, $w_{(1)} \in W_1$ and $w_{(2)} \in W_2,$ we
have the lower truncation condition
\[
(w_{(1)})_{n}w_{(2)} = 0\;\;  \mbox{for}\;\; n \;\; \mbox{whose real part is
sufficiently large;}
\]
the ``commutator formula''
\begin{eqnarray}\label{iocomm}
\lefteqn{Y_3(v,x_1){\cal Y}(w_{(1)},x_2)w_{(2)}-
{\cal Y}(w_{(1)},x_2)Y_2(v,x_1)w_{(2)}}\nno \\
&&{= \Res_{x_0}x^{-1}_2\delta \left( {x_1-x_0\over x_2}\right) {\cal Y}
(Y_1(v,x_0)w_{(1)},x_2) w_{(2)}};
\end{eqnarray}
and the $L(-1)$-derivative property
\begin{equation}\label{ioL(-1)}
{d\over dx}{\cal Y}(w_{(1)},x) = {\cal Y}(L(-1)w_{(1)},x),
\end{equation}
where  $L(-1)$  is the operator acting on  $W_{1}$. }
\end{defi}

\begin{rema}{\rm
For the notions of vertex (operator) algebra and module for a vertex
(operator) algebra, the term ``commutator formula'' has the intuitive
meaning---it is a formula for the commutator of two operators (acting
on the same space).  In the context of intertwining operators, even
though the similar formula, (\ref{iocomm}), does not involve a
commutator of two operators acting on the same space, we still call it
the ``commutator formula.''}
\end{rema}

Clearly, a quasi-intertwining operator ${\cal Y}$ is an intertwining
operator if and only if it further satisfies the Jacobi identity
\begin{eqnarray}\label{iojac}
\lefteqn{x^{-1}_0\delta \left( {x_1-x_2\over x_0}\right)
Y_3(v,x_1){\cal Y}(w_{(1)},x_2)w_{(2)}}\nno\\
&&\hspace{2em}- x^{-1}_0\delta \left( {x_2-x_1\over -x_0}\right)
{\cal Y}(w_{(1)},x_2)Y_2(v,x_1)w_{(2)}\nno \\
&&{= x^{-1}_2\delta \left( {x_1-x_0\over x_2}\right) {\cal Y}
(Y_1(v,x_0)w_{(1)},x_2) w_{(2)}}
\end{eqnarray}
for $v\in V$, $w_{(1)}\in W_1$ and $w_{(2)}\in W_2$; (\ref{iocomm}) of
course follows from (\ref{iojac}) by taking $\Res_{x_0}$.  It is clear
that the quasi-intertwining operators of the same type form a vector
space containing the space of intertwining operators as a subspace.

Recall from \cite{Mi} (see also \cite{HLZ1}, \cite{HLZ2}) the notion
of logarithmic intertwining operator:
\begin{defi}\label{logio}{\rm
Let $W_1$, $W_2$, $W_3$ be weak modules for a vertex operator algebra
$V$. A {\em logarithmic intertwining operator of type ${W_3\choose
W_1\,W_2}$} is a linear map
\[
{\cal Y}(\cdot, x)\cdot: W_1\otimes W_2\to W_3\{x\}[\log x],
\]
or equivalently,
\[
w_{(1)}\otimes w_{(2)}\mapsto{\cal Y}(w_{(1)},x)w_{(2)}=\sum_{n\in {\mathbb
C}}\sum_{k\in {\mathbb N}}{w_{(1)}}_{n;\,k}^{\cal Y}w_{(2)}x^{-n-1}(\log
x)^k\in W_3\{x\}[\log x]
\]
for all $w_{(1)}\in W_1$ and $w_{(2)}\in W_2$, such that the following
conditions are satisfied: the lower truncation condition: for
any $w_{(1)}\in W_1$, $w_{(2)}\in W_2$ and $k\in {\mathbb N}$,
\[
{w_{(1)}}_{n;\,k}^{\cal Y}w_{(2)}=0\;\;\mbox{ for }\;n
\;\mbox{ whose real part is sufficiently large;}
\]
the Jacobi identity (\ref{iojac}) for $v\in V$, $w_{(1)}\in W_1$
and $w_{(2)}\in W_2$ ; and the $L(-1)$-derivative property
(\ref{ioL(-1)}) for any $w_{(1)}\in W_1$. }
\end{defi}

\begin{rema}{\rm
The notion of logarithmic intertwining operator defined in \cite{HLZ1}
and \cite{HLZ2} is slightly more general than this one.  In this
paper, for brevity we adopt the original definition from \cite{Mi}
instead.}
\end{rema}

By analogy with the notion of quasi-intertwining operator, we can define
the notion of logarithmic quasi-intertwining operator, as follows:
\begin{defi}\label{logquasi}
{\rm
Let $W_1$, $W_2$, $W_3$ be weak modules for a vertex operator algebra
$V$. A {\em logarithmic quasi-intertwining operator of type ${W_3\choose
W_1\,W_2}$} is a linear map
\[
{\cal Y}(\cdot, x)\cdot: W_1\otimes W_2\to W_3\{x\}[\log x]
\]
that satisfies all the axioms in the definition of
logarithmic intertwining operator in Definition \ref{logio} except
that the Jacobi identity is replaced by the commutator formula
(\ref{iocomm}).  }
\end{defi}

{}From now on, unless otherwise stated, $(W_1,Y_1)$, $(W_2,Y_2)$
and $(W_3,Y_3)$ are assumed to be {\it
generalized} $V$-{\it modules} in the sense of \cite{HLZ1} and
\cite{HLZ2}, that is, weak $V$-modules satisfying all the axioms in
the definition of the notion of $V$-module (see \cite{FHL},
\cite{FLM}, \cite{LL} or \cite{tensor1}) except that the underlying
vector spaces are allowed to be direct sums of generalized
eigenspaces, not just eigenspaces, of the operator $L(0)$; in
particular, the $L(0)$-generalized eigenspaces are finite dimensional.
We refer the reader to \cite{HLZ1} and \cite{HLZ2} for basic notions
related to generalized $V$-modules. In particular, we have the notions
of algebraic completion and of contragredient module for a generalized
$V$-module.

In parallel to the notion of quasi-intertwining operator, we have:
\begin{defi}\label{qim}{\rm
A {\em quasi-$P(z)$-intertwining map of type ${W_{3}\choose
W_{1}W_{2}}$} is a linear map $F: W_{1}\otimes W_{2} \to
\overline{W}_{3}$ (the algebraic completion of $W_3$ with respect to
the grading by weights) satisfying the condition
\begin{eqnarray}\label{imcomm}
\lefteqn{Y_3(v,x_1)F(w_{(1)}\otimes w_{(2)})-F(w_{(1)}\otimes
Y_2(v,x_1)w_{(2)})}\nno\\
&&\hspace{3em}=\Res_{x_0}z^{-1}\delta\left(\frac{x_{1}-x_{0}}{z}\right)
F(Y_1(v,x_0)w_{(1)}\otimes w_{(2)})
\end{eqnarray}
for $v\in V$, $w_{(1)}\in W_{1}$, $w_{(2)}\in W_{2}$.
}
\end{defi}

Note that the left-hand side of (\ref{imcomm}) is well defined, by the
same argument as was used for the left-hand side of formula (4.2) in
\cite{tensor1}; that argument indeed remains valid for generalized
modules.

Clearly, a quasi-$P(z)$-intertwining map ${\cal Y}$ of type
${W_{3}\choose W_{1}\,W_{2}}$ is a $P(z)$-intertwining map if and only
if it further satisfies the Jacobi identity
\begin{eqnarray}\label{imjac}
\lefteqn{x_{0}^{-1}\delta\left(\frac{ x_{1}-z}{x_{0}}\right)
Y_{3}(v, x_{1})F(w_{(1)}\otimes w_{(2)})=}\nonumber\\
&&=z^{-1}\delta\left(\frac{x_{1}-x_{0}}{z}\right)
F(Y_{1}(v, x_{0})w_{(1)}\otimes w_{(2)})\nonumber\\
&&\hspace{2em}+x_{0}^{-1}\delta\left(\frac{z-x_{1}}{-x_{0}}\right)
F(w_{(1)}\otimes Y_{2}(v, x_{1})w_{(2)})
\end{eqnarray}
for $v\in V$, $w_{(1)}\in W_{1}$, $w_{(2)}\in W_{2}$; (\ref{imcomm})
follows from (\ref{imjac}) by taking $\Res_{x_0}$.
(It is important to keep in mind that the left-hand side of
(\ref{imjac}) is well defined, by the considerations at the beginning
of Section 4 of \cite{tensor1}.  Clearly, the
quasi-$P(z)$-intertwining maps of the same type form a vector space
containing the space of $P(z)$-intertwining maps as a subspace.

In case $W_1$, $W_2$ and $W_3$ are ordinary $V$-modules, given a fixed
integer $p$, by analogy with the maps defined in (12.3) and (12.4) in
\cite{tensor3}, we have the following maps between the spaces of
quasi-intertwining operators and of quasi-$P(z)$-intertwining maps of
the same type: For a quasi-intertwining operator ${\cal Y}$ of type
${W_3\choose W_1\, W_2}$, define $F_{{\cal Y},p}: W_1\otimes
W_2\to\overline{W}_3$ by
\[
F_{{\cal Y},p}(w_{(1)}\otimes w_{(2)})={\cal Y}(w_{(1)},e^{l_p(z)})
w_{(2)}
\]
for all $w_{(1)}\in W_1$ and $w_{(2)}\in W_2$, where we follow the
notation
\begin{eqnarray*}
\log z&=&\log |z| +i \, {\rm arg\,}
z\;\;\mbox{for}\;0\leq {\rm arg\,}z<2\pi,\\
l_p(z)&=&\log z+2\pi ip, \ p\in\Z
\end{eqnarray*}
in \cite{tensor1} for branches of the $\log$ function.  On the other
hand, let $F$ be a quasi-$P(z)$-intertwining map of type ${W_3\choose
W_1\, W_2}$.  For homogeneous elements $w_{(1)}\in W_1$ and
$w_{(2)}\in W_2$ and $n\in\C$, define $(w_{(1)})_n w_{(2)}$ to be the
projection of $F(w_{(1)}\otimes w_{(2)})$ to the homogeneous subspace
of $W_3$ of weight $\wt w_{(1)}-n-1+\wt w_{(2)}$ multiplied by
$e^{(n+1)l_p(z)}$, and define
\[
{\cal Y}_{F,p}(w_{(1)},x)w_{(2)}=\sum_{n\in\C}(w_{(1)})_n
w_{(2)}x^{-n-1};
\]
then extend by linearity to define ${\cal Y}_{F,p}: W_1\otimes W_2\to
{W}_3 \{x\}$.

It was shown in Proposition 12.2 of \cite{tensor3} (see also
Proposition 4.7 in \cite{tensor1}) that these two maps give linear
isomorphisms between the space of intertwining operators and the space
of $P(z)$-intertwining maps of the same type. By replacing all Jacobi
identities by the corresponding commutator formulas in the proof, we
see that these two maps also give linear isomorphisms between the
space of quasi-intertwining operators and the space of
quasi-$P(z)$-intertwining maps of the same type.  (The straightforward
argument is carried out in \cite{HLZ2}.)  That is, we have:

\begin{propo}\label{ioim}
Assume that $W_1$, $W_2$ and $W_3$ are ordinary $V$-modules.
For $p\in\Z$, the correspondence ${\cal Y}\mapsto F_{{\cal Y},p}$ is a
linear isomorphism {from} the space of quasi-intertwining operators of
type ${W_{3}}\choose {W_{1}\; W_{2}}$ to the space of
quasi-$P(z)$-intertwining maps of the same type.  Its inverse map is
given by $F\mapsto {\cal Y}_{F, p}$.  \epf
\end{propo}

More generally, if $W_1$, $W_2$ and $W_3$ are generalized (rather than
ordinary) $V$-modules in the sense of \cite{HLZ1} and \cite{HLZ2},
then following the argument in \cite{HLZ2} we have a result similar to
Proposition \ref{ioim} giving the correspondence between the
quasi-$P(z)$-intertwining maps and the logarithmic quasi-intertwining
operators.

Here is an easy example of a quasi-$P(z)$-intertwining map that is not
a $P(z)$-intertwining map:
\begin{exam}\label{badim}{\rm
Take $V$ to be the vertex operator algebra constructed from a
finite-dimen\-sional commutative associative algebra with identity, with
the vertex operator map defined by $Y(a,x)b=ab$ for $a,b\in V$, with
the vacuum vector ${\bf 1}$ taken to be $1$ and with $\omega = 0$.
Then since
the notion of module for a vertex algebra can be characterized in
terms of an associativity property, the modules for $V$ as a vertex
operator algebra are precisely the finite-dimensional modules for $V$
as an associative algebra (see \cite{B}; cf. \cite{LL}). For two
$V$-modules $W_1$ and $W_2$, the vector space $W_1\otimes W_2$ is a
$V$-module under the action given by
\begin{equation}\label{ca-y}
Y(v,x)(w_{(1)}\otimes w_{(2)})\;(=v\cdot(w_{(1)}\otimes w_{(2)}))=
w_{(1)}\otimes(v\cdot w_{(2)})
\end{equation}
for $v\in V$, $w_{(1)}\in W_{1}$ and $w_{(2)}\in W_{2}$.
The identity map on $W_1\otimes W_2$ is a
quasi-$P(z)$-intertwining map of type ${W_1\otimes W_2\choose
W_1\;W_2}$ (formula (\ref{imcomm}) being just (\ref{ca-y})
itself). However it is not a $P(z)$-intertwining map because the
Jacobi identity demands that
\[
(v\cdot w_{(1)})\otimes w_{(2)}=w_{(1)}\otimes(v\cdot w_{(2)})
\]
for any $v\in V$, $w_{(1)}\in W_1$ and $w_{(2)}\in W_2$, which of
course is not true in general.  See Remark \ref{badlambda} below for a
further discussion of this (counter)example.}
\end{exam}

\br{rbadio} {\rm By the above, this example of course immediately
gives a quasi-intertwining operator that is not an intertwining
operator.}  \er

In the following we will sometimes use results from \cite{HLZ1} and
\cite{HLZ2} for our generalized $V$-modules $W_1$, $W_2$ and $W_3$,
but the reader may simply take $W_1$, $W_2$ and $W_3$ to be ordinary
$V$-modules.

Just as in \cite{tensor1}, formulas (3.4) and (3.5), set
$$Y_{t}(v, x)=v\otimes x^{-1}\delta\left(\frac{t}{x}\right)$$
for $v\in V$.  Also, just as in \cite{tensor1}, formula (3.20),
for a generalized $V$-module $(W, Y)$, write
$$Y^{o}(v, x)=Y(e^{xL(1)}(-x^{-2})^{L(0)}v, x^{-1})$$
for $v\in V$.
(Here we are using the notation $Y^{o}$, as in \cite{HLZ1} and
\cite{HLZ2}, rather than the original notation $Y^*$ of
\cite{tensor1}.)
Also, let $\iota_{+}: \C[t,t^{- 1}, (z^{-1}-t)^{-1}]
\to \C((t))$ (the space of formal Laurent series in $t$
with only finitely many negative powers of $t$) be the natural map.
As in formula (13.2) of \cite{tensor3},
we define a linear action $\tau_{P(z)}$ of the space
\begin{equation}\label{vspace}
V \otimes \iota_{+}\C[t,t^{- 1}, (z^{-1}-t)^{-1}]
\end{equation}
on $(W_1\otimes W_2)^*$ by
\begin{eqnarray}\label{tau0}
\lefteqn{\left(\tau_{P(z)}
\left(x_{0}^{-1}\delta\left(\frac{x^{-1}_{1}-z}{x_{0}}\right)
Y_{t}(v, x_{1})\right)\lambda\right)(w_{(1)}\otimes w_{(2)})=}\nonumber\\
&&=z^{-1}\delta\left(\frac{x^{-1}_{1}-x_{0}}{z}\right)
\lambda(Y_{1}(e^{x_{1}L(1)}(-x_{1}^{-2})^{L(0)}v, x_{0})w_{(1)}\otimes w_{(2)})
\nonumber\\
&&\quad +x^{-1}_{0}\delta\left(\frac{z-x^{-1}_{1}}{-x_{0}}\right)
\lambda(w_{(1)}\otimes Y_{2}^{o}(v, x_{1})w_{(2)})
\end{eqnarray}
for $v\in V$, $\lambda\in (W_{1}\otimes W_{2})^{*}$, $w_{(1)}\in
W_{1}$ and $w_{(2)}\in W_{2}$.  As in \cite{tensor3}, this formula
does indeed give a well-defined linear action (in generating-function
form) of the space (\ref{vspace}); see Section 3 of \cite{tensor1}.
This action (\ref{tau0}) restricts in particular to an
action of $V \otimes \C[t,t^{-1}]$ on $(W_1\otimes
W_2)^*$, given in generating-function form by
$\tau_{P(z)}(Y_{t}(v,x))$; one takes the residue with respect to $x_0$
of both sides of (\ref{tau0}).

We will write $W_1\warning W_2$, which can be read ``$W_1$ warning
$W_2$'', for the largest weak $V$-module inside $(W_1\otimes W_2)^*$
with respect to the action $\tau_{P(z)}$ of $V \otimes
\C[t,t^{-1}]$. (Here we omit the information about $P(z)$ {}from the
notation.) It is clear that $W_1\warning W_2$ does exist and equals
the sum (or union) of all weak $V$-modules inside $(W_1\otimes
W_2)^*$.  Of course, all the elements of $W_1\warning W_2$ satisfy the
lower truncation condition and the Jacobi identity with respect to the
action $\tau_{P(z)}$.  (Warning: This space $W_1\warning W_2$ can be
strictly larger than the subspace $W_1\hboxtr_{P(z)} W_2$ of
$(W_1\otimes W_2)^*$ defined in formula (13.13) of \cite{tensor3}, as
we will show below.)

Denote by $W'$ the contragredient module of a generalized $V$-module
$W$, defined by exactly the same procedure as was carried out in
\cite{FHL} for ordinary (as opposed to generalized) modules.
For a linear map $F$ from
$W_1\otimes W_2$ to $\overline{W}_3$, define a linear map
$F^\vee:W'_3\to(W_1\otimes W_2)^*$ by
\begin{equation}\label{FV}
\langle F^\vee(\alpha),w_{(1)}\otimes w_{(2)}\rangle_{W_1\otimes W_2}=
\langle\alpha,F(w_{(1)}\otimes w_{(2)})\rangle_{\overline{W}_3}
\end{equation}
for $\alpha\in W'_3$, $w_{(1)}\in W_1$ and $w_{(2)}\in W_2$.
(The subscripts of course designate the pairings; sometimes we will
omit such subscripts.)  For the
case in which $W_1$, $W_2$ and $W_3$ are ordinary modules, it was
observed in \cite{tensor3}, Proposition 13.1 (see also \cite{tensor1},
Proposition 5.3) that $F$ is a $P(z)$-intertwining map of type
${W_{3}}\choose {W_{1}\; W_{2}}$ if and only if $F^\vee$ intertwines
the two actions of the space (\ref{vspace}) on $W'_3$ (on which a
monomial $v\otimes t^n$ acts as $v_n$ and a general element acts
according to the action of each of its monomials) and on $(W_1\otimes
W_2)^*$ (by formula (\ref{tau0})). The same observation still holds for the
case of generalized modules (cf. \cite{HLZ2}).  For
quasi-$P(z)$-intertwining maps, we shall prove:

\begin{propo}\label{prop-corr}
The map $F$ is a quasi-$P(z)$-intertwining map of
type ${W_{3}}\choose {W_{1}\;
W_{2}}$ if and only if $F^\vee$ intertwines the two actions of
$V\otimes\C[t,t^{-1}]$ on $W'_3$ and $(W_1\otimes W_2)^*$.
\end{propo}

\begin{rema}{\rm
In the statement of Proposition \ref{prop-corr} we have avoided saying that
$F^\vee$ is a $V$-homomorphism since the target space,
$(W_1\otimes W_2)^*$, is rarely a (generalized) $V$-module.  }
\end{rema}

Before giving the proof, we first write the action (\ref{tau0}) in an
alternative form, which will be more convenient in this paper, as
follows:
\begin{eqnarray}\label{tau}
\lefteqn{\langle\tau_{P(z)}
\left(x_{0}^{-1}\delta\left(\frac{x_{1}-z}{x_{0}}\right) Y^o_{t}(v,
x_{1})\right)\lambda, w_{(1)}\otimes w_{(2)}\rangle=}\nno\\
&&=\langle\lambda,z^{-1}\delta\left(\frac{x_{1}-x_{0}}{z}\right)Y_{1}(v,
x_{0})w_{(1)}\otimes w_{(2)} \nno\\ &&\qquad
+x^{-1}_{0}\delta\left(\frac{z-x_{1}}{-x_{0}}\right) w_{(1)}\otimes
Y_{2}(v, x_{1})w_{(2)}\rangle
\end{eqnarray}
for $v\in V$, $\lambda\in (W_{1}\otimes W_{2})^{*}$, $w_{(1)}\in
W_{1}$ and $w_{(2)}\in W_{2}$, where $Y^o_{t}(v,x_{1})$ is defined
by
\[
Y_t^o(v,x)=Y_t(e^{xL(1)}(-x^{-2})^{L(0)}v,x^{-1})
=e^{xL(1)}(-x^{-2})^{L(0)}v\otimes x
\delta\left(\frac{t}{x^{-1}}\right),
\]
which in turn equals
$$(-1)^{\swt v}\sum_{m\in\N}\frac 1{m!}(L(1)^mv)\otimes
t^{-m-2+2\swt v} x^{-1}\delta\bigg(\frac {t^{-1}}{x}\bigg)$$
in case $v$ is homogeneous, by formulas (3.25), (3.30), (3.32) and
(3.38) of \cite{tensor1}.
The equivalence of (\ref{tau0}) and (\ref{tau}) can be seen by first
replacing $x_1$ by $x_1^{-1}$ and then replacing $v$ by
$e^{x_1L(1)}(-x_1^{-2})^{L(0)}v$ in either direction (recall
Proposition 5.3.1 of \cite{FHL}).

\medskip\noindent {\it Proof of Proposition \ref{prop-corr}\quad} For
any linear map $F$ from $W_1\otimes W_2$ to $\overline{W}_3$, the
condition that the map
$F^\vee$ defined by (\ref{FV}) intertwines the two actions of
$V\otimes\C[t,t^{-1}]$ on
$W'_3$ and $(W_1\otimes W_2)^*$ means exactly that for any $\alpha\in
W'_3$, $v\in V$, $w_{(1)}\in W_1$ and $w_{(2)}\in W_2$,
\begin{equation}\label{intw}
\langle F^\vee((Y'_3)^o(v,x_1)\alpha),w_{(1)}\otimes w_{(2)}
\rangle_{W_1\otimes W_2}
=\langle\tau_{P(z)}(Y^o_{t}(v,x_{1}))
F^\vee(\alpha),w_{(1)}\otimes w_{(2)}\rangle_{W_1\otimes W_2}.
\end{equation}
The left-hand side of (\ref{intw}) is
\begin{eqnarray*}
\langle F^\vee((Y'_3)^o(v,x_1)\alpha),w_{(1)}\otimes w_{(2)}
\rangle_{W_1\otimes W_2}&=&
\langle (Y'_3)^o(v,x_1)\alpha,F(w_{(1)}\otimes w_{(2)})
\rangle_{\overline{W}_3}\\
&=&\langle\alpha,Y_3(v,x_1)F(w_{(1)}\otimes w_{(2)})
\rangle_{\overline{W}_3},
\end{eqnarray*}
while by setting $\lambda = F^\vee(\alpha)$ in (\ref{tau}) and
then taking $\Res_{x_0}$, we see that the right-hand side of
(\ref{intw}) is
\begin{eqnarray*}
\lefteqn{\langle\tau_{P(z)}(Y^o_{t}(v,x_{1}))F^\vee(\alpha),w_{(1)}
\otimes w_{(2)}\rangle_{W_1\otimes W_2}=}\\
&&\hspace{6em}=\langle F^\vee(\alpha),\Res_{x_0}z^{-1}\delta\left
(\frac{x_{1}-x_{0}}{z}\right)Y_{1}(v, x_{0})w_{(1)}\otimes w_{(2)}\\
&&\hspace{12em}+w_{(1)}\otimes Y_{2}(v, x_{1})w_{(2)}
\rangle_{W_1\otimes W_2}\\
&&\hspace{6em}=\langle\alpha,\Res_{x_0}z^{-1}\delta\left(\frac{x_{1}-
x_{0}}{z}\right)F(Y_{1}(v,x_{0})w_{(1)}\otimes w_{(2)})\\
&&\hspace{12em}+F(w_{(1)}\otimes Y_{2}(v,x_{1})w_{(2)})
\rangle_{\overline{W}_3}.
\end{eqnarray*}
The proposition follows immediately. \epf

\br{Jacobionlambda}
{\rm Note that for fixed $\lambda$, $\tau_{P(z)}(Y_{t}(v,x))\lambda$ is
lower truncated (with respect to $x$) for any $v\in V$ if and only if
$\tau_{P(z)}(Y^o_{t}(v,x))\lambda$ is upper truncated for any $v\in V$.
Moreover, in this case, the Jacobi identity for $\tau_{P(z)}(Y_{t}(\cdot,x))$
holds on $\lambda$ if and only if the opposite Jacobi identity for
$\tau_{P(z)}(Y^o_{t}(\cdot,x))$ (see formula (3.23) in \cite{tensor1})
holds on $\lambda$.  Indeed, first assume that the Jacobi identity for
$\tau_{P(z)}(Y_{t}(\cdot,x))$ holds on $\lambda$.  An examination of
the proof of Theorem 5.2.1 in \cite{FHL}, which asserts that the
contragredient of a module is indeed a module, in fact proves the
desired opposite Jacobi identity.  (A similar observation was made in
reference to formula (3.23) in \cite{tensor1}.)  For the converse, one
sees that the relevant steps in the proof of Theorem 5.2.1 in
\cite{FHL} are reversible.}
\er

\begin{thm}\label{Fv(W)}
Let $W_1$, $W_2$ be generalized $V$-modules and $W_3$ be an ordinary
(respectively, generalized) $V$-module. Let $F$ be a
quasi-$P(z)$-intertwining map of type ${W_{3}}\choose
{W_{1}\; W_{2}}$. Then for any $\alpha\in W'_3$,
$F^\vee(\alpha)\in(W_1\otimes W_2)^*$ satisfies the lower truncation
condition and the Jacobi identity with respect to the action
$\tau_{P(z)}$. In particular, $F^\vee(W'_3)\subset(W_1\otimes W_2)^*$
is an ordinary (respectively, generalized) $V$-module
and $F^\vee : W'_3 \rightarrow F^\vee (W'_3)$ is a module map
(respectively, a map of generalized modules).  Conversely,
every ordinary (respectively, generalized) $V$-module inside
$(W_1\otimes W_2)^*$ arises in this way.
\end{thm}
\pf Let $F$ be as in the assumption.  Then for any $\alpha\in W'_3$, by
Proposition \ref{prop-corr} we have
\begin{equation}\label{intw2}
\tau_{P(z)}(Y^o_{t}(v,x_{1}))F^\vee(\alpha)=F^\vee((Y'_3)^{o}(v,x_1)\alpha)
\end{equation}
for any $v\in V$. Since the right-hand side of (\ref{intw2}) is upper
truncated in
$x_1$, so is the left-hand side. Hence we have the lower truncation
condition with respect to the action $\tau_{P(z)}$. The Jacobi
identity on $F^\vee(\alpha)$
follows from (\ref{intw2}) and the fact that $\alpha$
satisfies the Jacobi identity on $W'_3$
(recall Remark \ref{Jacobionlambda}).
Also, $\tau_{P(z)}(Y_{t}({\bf 1},x))=1$ from the definitions.
Therefore, $F^\vee(W'_3)$ is a
weak $V$-module. But as an image of the ordinary (respectively,
generalized) $V$-module $W'_3$, it must be an ordinary (respectively,
generalized) $V$-module itself.

Conversely, let $M$ be a subspace of $(W_1\otimes W_2)^*$ that becomes
an ordinary (respectively, generalized) $V$-module under the action
$\tau_{P(z)}$ of $V\otimes\C[t,t^{-1}]$. Take $W_3=M'$, the
contragredient module of $M$, and define $F:W_1\otimes W_2\to
\overline{W}_3$ by
\begin{equation}\label{getf}
\langle\alpha,F(w_{(1)}\otimes w_{(2)})\rangle_{\overline{W}_3}=
\langle\alpha,w_{(1)}\otimes w_{(2)}\rangle_{W_1\otimes W_2}
\end{equation}
for any $\alpha\in W'_3=M\subset (W_1\otimes W_2)^*$, $w_{(1)}\in
W_1$ and $w_{(2)}\in W_2$. By using $\Res_{x_0}$ of (\ref{tau})
we have
\begin{eqnarray*}
\lefteqn{\langle\alpha,Y_3(v,x_1)F(w_{(1)}\otimes w_{(2)})
\rangle_{\overline{W}_3}}\nno\\
&&={\langle Y'_3(v,x_1)\alpha,F(w_{(1)}\otimes w_{(2)})
\rangle_{\overline{W}_3}}\nno\\
&&=\langle\tau_{P(z)}(Y^o_{t}(v,x_{1}))\alpha,F(w_{(1)} \otimes
w_{(2)})\rangle_{\overline{W}_3}\nno\\
&&=\langle\tau_{P(z)}(Y^o_{t}(v,x_{1}))\alpha,w_{(1)} \otimes
w_{(2)}\rangle_{W_1\otimes W_2}\nno\\
&&=\langle\alpha,\Res_{x_0}z^{-1}\delta\left
(\frac{x_{1}-x_{0}}{z}\right)Y_{1}(v, x_{0})w_{(1)}\otimes
w_{(2)}+w_{(1)}\otimes Y_{2}(v, x_{1})w_{(2)}\rangle_{W_1\otimes
W_2}\nno\\
&&=\langle\alpha,\Res_{x_0}z^{-1}\delta\left(\frac{x_{1}-x_{0}}{z}
\right)F(Y_{1}(v,x_{0})w_{(1)}\otimes w_{(2)})+F(w_{(1)}\otimes
Y_{2}(v,x_{1})w_{(2)})\rangle_{\overline{W}_3}\qquad\nno
\end{eqnarray*}
for any $\alpha\in M$, $v\in V$, $w_{(1)}\in W_1$ and $w_{(2)}\in W_2$. This
shows that $F$ is a quasi-$P(z)$-intertwining map. In addition, by
(\ref{getf}) we also have that
\[
\langle F^\vee(\alpha),w_{(1)}\otimes w_{(2)}\rangle_{W_1\otimes W_2}=
\langle \alpha,F(w_{(1)}\otimes w_{(2)})\rangle_{\overline{W}_3}=
\langle \alpha,w_{(1)}\otimes w_{(2)}\rangle_{W_1\otimes W_2}
\]
for any $\alpha\in W'_3=M$, $w_{(1)}\in W_1$ and $w_{(2)}\in
W_2$, so that $F^\vee$ is the identity map on $W'_3=M$, and
$M=F^\vee(W'_3)$.
\epfv

An immediate consequence is:
\begin{corol}\label{cor1}
Suppose $W_4$ is also a generalized $V$-module (in addition to $W_1$,
$W_2$ and $W_3$).  Let
$F_1$ and $F_2$ be quasi-$P(z)$-intertwining maps of types
${W_{3}}\choose {W_{1}\; W_{2}}$ and ${W_{4}}\choose {W_{1}\; W_{2}}$,
respectively. Assume that both $W_3$ and $W_4$ are irreducible. Then
the generalized $V$-modules $F_1^\vee(W'_3)$ and $F_2^\vee(W'_4)$ inside
$(W_1\otimes W_2)^*$ are irreducible and in particular
either coincide with each other or have
intersection $0$. \epf
\end{corol}

\begin{rema}\label{commonker}{\rm
In the first half of Theorem \ref{Fv(W)}, let $W\subset W_3$ be the space
spanned by the homogeneous components of the elements
$F(w_{(1)}\otimes w_{(2)})\in\overline{W}_3$ for
all $w_{(1)}\in W_1$ and $w_{(2)}\in W_2$. Then by formula
(\ref{imcomm}) it is clear that $W$ is closed under the action of each
$v_n$ for $v\in V$ and $n\in\Z$. Hence, in
the case that $W_3$ is an ordinary (respectively, generalized) $V$-module, $W$
is an ordinary (respectively, generalized) $V$-submodule of
$W_3$. Furthermore, $F^\vee(W'_3)$ is naturally isomorphic
to $W'$ as an ordinary (respectively, generalized) $V$-module; indeed,
both $V$-homomorphisms
$W'_3\to F^\vee(W'_3)$ and $W'_3\to W'$ are surjective and have the
common kernel
\[
\{\alpha\in W'_3\,|\,\langle \alpha,W\rangle_{W_3}=0\}
\]
(recall (\ref{FV})).
}
\end{rema}

\begin{rema}{\rm
Consider the special case of Remark \ref{commonker} in which $W_1$ and
$W_2$ are ordinary modules, $W_3$ is $W_1\boxtimes_{P(z)}
W_2=(W_1\hboxtr_{P(z)}W_2)'$ if it exists, and $F$ is the canonical
$P(z)$-intertwining map $W_1\otimes
W_2\to\overline{W_1\boxtimes_{P(z)} W_2}$ coming from the canonical
injection $F^\vee: W_1\hboxtr_{P(z)} W_2\to(W_1\otimes W_2)^*$. Then
the map $W'_3\to F^\vee(W'_3)=W'_3$ is the identity and so the map
$W'_3\to W'$ is an isomorphism of modules. In particular, $W=W_3$, and
we have recovered Lemma 14.9 of \cite{tensor4}: The homogeneous
components of the tensor product elements $w_{(1)}\boxtimes_{P(z)}
w_{(2)}\; (=F(w_{(1)}\otimes w_{(2)}))$ span the tensor product module
$W_1\boxtimes_{P(z)} W_2$.  }
\end{rema}

Recall from \cite{tensor3}, formulas (13.12) and (13.16),
that an element $\lambda\in(W_1\otimes
W_2)^*$ is said to satisfy the {\em $P(z)$-compatibility condition} if
$\lambda$ satisfies the lower truncation condition with respect to the
action $\tau_{P(z)}$ and for any $v\in V$, the following formula
holds:
\begin{equation}\label{comp0}
\tau_{P(z)}\left(x_{0}^{-1}\delta\left(\frac{x^{-1}_{1}-z}{x_{0}}
\right) Y_{t}(v, x_{1})\right)\lambda=
x_{0}^{-1}\delta\left(\frac{x^{-1}_{1}-z}{x_{0}}\right)
\tau_{P(z)}(Y_{t}(v, x_{1}))\lambda;
\end{equation}
that is, the action of the space (\ref{vspace}) on $\lambda$ given by
(\ref{tau0}) is compatible with the action of the space
$V\otimes\C[t,t^{-1}]$ on $\lambda$ given by restricting (\ref{tau0})
to the elements $Y_{t}(v, x_{1})$.
Recall also that a subspace of $(W_1\otimes W_2)^*$
(in particular, a generalized $V$-module inside $(W_1\otimes W_2)^*$)
is said to
be {\em $P(z)$-compatible} if all of its elements satisfy the
$P(z)$-compatibility condition. By Theorem \ref{Fv(W)} we have:

\begin{corol}\label{cor2}
Let $F$ be a quasi-$P(z)$-intertwining map of type ${W_{3}}\choose
{W_{1}\; W_{2}}$. Assume that $W_3$ is irreducible. Then if
$F^\vee(W'_3)$ is not $P(z)$-compatible, none of its nonzero elements
satisfies the $P(z)$-compatibility condition.
\end{corol}
\pf If $F^\vee(W'_3)=0$ the conclusion is trivial. Otherwise, by
Theorem \ref{Fv(W)} and the irreducibility of $W'_3$ we see that
$F^\vee(W'_3)$ is an irreducible generalized $V$-module. The statement
now follows {}from the fact that the set of elements satisfying the
$P(z)$-compatibility condition is stable under the action
$\tau_{P(z)}$ (see Theorem 13.9 of \cite{tensor3}, Proposition
6.2 of \cite{tensor1}, and their generalization in \cite{HLZ2} for the
case of generalized modules). \epfv

We have:
\begin{thm}\label{pcorresp}
Let $F$ be a quasi-$P(z)$-intertwining map of type ${W_{3}}\choose
{W_{1}\; W_{2}}$.  Then the generalized module (ordinary if ${W_{3}}$
is ordinary)
$F^\vee(W'_3)$ is $P(z)$-compatible if and only
if $F$ is in fact a $P(z)$-intertwining map.
\end{thm}
\pf For convenience we use an equivalent form of (\ref{comp0}) as
follows:
\begin{equation}\label{comp}
\tau_{P(z)} \left(x_{0}^{-1}\delta\left(\frac{x_{1}-z}{x_{0}}\right)
Y^o_{t}(v, x_{1})\right)\lambda=x_{0}^{-1}\delta
\left(\frac{x_{1}-z}{x_{0}}\right)\tau_{P(z)}(Y^o_{t}(v,x_{1}))
\lambda
\end{equation}
(recall the equivalence between (\ref{tau0}) and (\ref{tau})).
Let $F$ be a quasi-$P(z)$-intertwining map. In (\ref{comp}), setting
$\lambda=F^\vee(\alpha)$ and applying to $w_{(1)}\otimes w_{(2)}$ for
$\alpha\in W'_3$, $w_{(1)}\in W_1$ and $w_{(2)}\in W_2$, we see that
the left-hand side becomes
\begin{eqnarray}\label{lcomp}
\lefteqn{\langle\tau_{P(z)}\left(x_{0}^{-1}\delta\left(\frac{x_{1}-z}
{x_{0}}\right) Y^o_{t}(v,x_{1})\right)F^\vee(\alpha),
w_{(1)}\otimes w_{(2)}\rangle=}\nno\\
&&\hspace{6em}=\langle F^\vee(\alpha),z^{-1}\delta\left(\frac{x_{1}
-x_{0}}{z}\right)Y_{1}(v, x_{0})w_{(1)}\otimes w_{(2)}\nno\\
&&\hspace{12em}+x^{-1}_{0}\delta\left(\frac{z-x_{1}}{-x_{0}}\right)
w_{(1)}\otimes Y_{2}(v, x_{1})w_{(2)}\rangle\nno\\
&&\hspace{6em}=\langle\alpha,z^{-1}\delta\left(\frac{x_{1}-x_{0}}{z}
\right)F(Y_{1}(v,x_{0})w_{(1)}\otimes w_{(2)})\nno\\
&&\hspace{12em}+x^{-1}_{0}\delta\left(\frac{z-x_{1}}{-x_{0}}\right)
F(w_{(1)}\otimes Y_{2}(v, x_{1})w_{(2)})\rangle,
\end{eqnarray}
while the right-hand side becomes
\begin{eqnarray}\label{rcomp}
\lefteqn{x_{0}^{-1}\delta\left(\frac{x_{1}-z}{x_{0}}\right)\langle
\tau_{P(z)}(Y^o_{t}(v,x_{1}))F^\vee(\alpha),w_{(1)}
\otimes w_{(2)}\rangle=}\nno\\
&&\hspace{6em}=x_{0}^{-1}\delta\left(\frac{x_{1}-z}{x_{0}}\right)
\langle F^\vee(\alpha),\Res_{x_0}z^{-1}\delta\left(\frac{x_{1}
-x_{0}}{z}\right)Y_{1}(v, x_{0})w_{(1)}\otimes w_{(2)}\nno\\
&&\hspace{12em}+w_{(1)}\otimes Y_{2}(v, x_{1})w_{(2)}\rangle\nno\\
&&\hspace{6em}=x_{0}^{-1}\delta\left(\frac{x_{1}-z}{x_{0}}\right)
\langle\alpha,\Res_{x_0}z^{-1}\delta\left(\frac{x_{1}-x_{0}}{z}\right)
F(Y_{1}(v,x_{0})w_{(1)}\otimes w_{(2)})\nno\\
&&\hspace{12em}+F(w_{(1)}\otimes Y_{2}(v,x_{1})w_{(2)})\rangle\nno\\
&&\hspace{6em}=\langle\alpha,x_{0}^{-1}\delta\left(\frac{x_{1}-z}{x_{0}}
\right)Y_3(v,x_1)F(w_{(1)}\otimes w_{(2)})\rangle,
\end{eqnarray}
where in the last step, we have used (\ref{imcomm}).
Thus $F^\vee(W'_3)$ is $P(z)$-compatible if and only if for any
$\alpha\in W'_3$, the right-hand side of (\ref{lcomp}) is equal to the
right-hand side of (\ref{rcomp}) for any $w_{(1)}\in W_1$ and
$w_{(2)}\in W_2$, that is, if and only if (\ref{imjac}) is true for
any $w_{(1)}\in W_1$ and $w_{(2)}\in W_2$, which is equivalent to $F$
being a $P(z)$-intertwining map. \epf

\begin{rema}{\rm
Suppose that an element $\lambda$ of $(W_1\otimes W_2)^*$ satisfies the
$P(z)$-compatibility condition and generates a generalized $V$-module
$W_\lambda$ under the action $\tau_{P(z)}$ (cf.\ the $P(z)$-local
grading-restriction condition in \cite{tensor3}). Then $W_\lambda$ is
compatible, just as in the proof of Corollary \ref{cor2}. By Theorem
\ref{Fv(W)}, $W_\lambda=F^\vee(W')$ for some generalized $V$-module $W$
and quasi-$P(z)$-intertwining map $F$ of type ${W\choose
W_1\,W_2}$. Then Theorem \ref{pcorresp} ensures that $F$ is in fact a
$P(z)$-intertwining map. In particular, $\lambda$ lies in the image of
$F^\vee$ for the $P(z)$-intertwining map $F$.  }
\end{rema}

\begin{rema}\label{badlambda}{\rm
Theorem \ref{pcorresp} and Example \ref{badim} together provide
examples of non-$P(z)$-compatible modules.\footnote{These examples
show in particular that the
construction of the tensor product functor in math.QA/0309350---the
formula below formula (4.1)---appears to be incorrect:
As a consequence of the definition of $P(z)$-tensor
product (adopted from \cite{tensor1}), the contragredient of the
$P(z)$-tensor product of modules $W_1$ and $W_2$ is the union (or sum)
of all {\em $P(z)$-compatible} modules, rather than {\em all} modules,
inside $(W_1\otimes W_2)^*$. The examples we give here and below show
that the space defined in math.QA/0309350 is sometimes strictly larger than
the correct contragredient module of the tensor product module. In
particular, the arguments in math.QA/0309350 purport to
establish an assertion equivalent to associativity for
{\em quasi-intertwining} operators, which is not true. The correct
result, proved (in the logarithmic context) in \cite{HLZ1},
\cite{HLZ2}, is the associativity for {\it intertwining} operators;
this work generalizes the arguments in
\cite{tensor1}--\cite{tensor3}, \cite{tensor4} and of course is based
on the compatibility condition. For the examples in the present remark,
even when $W_1$ and $W_2$ are taken to be $V$ itself,
the construction in math.QA/0309350 results in a space strictly larger (in the
sense of homogeneous subspaces) than the correct tensor product,
$V$ itself. All of this illustrates why the compatibility condition
of \cite{tensor1}--\cite{tensor3} is crucial and cannot in general be
omitted. As we have mentioned, the compatibility condition remains
crucial in the
construction of the natural associativity isomorphisms among triple
tensor products in \cite{tensor4}, and the proofs of their fundamental
properties.  In the tensor product theory in
\cite{tensor1}--\cite{tensor3} and \cite{tensor4}, the compatibility
condition on elements of $(W_1\otimes W_2)^*$ was not a restriction on
the applicability of the theory; rather, it was a necessary condition
for obtaining the (correct) theory, and the same is certainly true for the
still more subtle logarithmic generalization of the theory in
\cite{HLZ1}--\cite{HLZ2}.}
It is instructive to write down the details of these
(counter)examples:  Take $V$ to be the vertex operator algebra constructed
{}from a finite-dimensional commutative associative algebra with identity
as in
Example \ref{badim}. For any $V$-modules $W_1$ and $W_2$ (that is,
finite-dimensional modules $W_1$ and $W_2$ for $V$ viewed as an
associative algebra), the action
$\tau_{P(z)}$ of $V\otimes \C[t,t^{-1}]$ on $(W_1\otimes W_2)^*$ is
given by:
\[
(\tau_{P(z)}(Y_{t}(v,x))\lambda)(w_{(1)}\otimes w_{(2)})=
\lambda(w_{(1)}\otimes(v\cdot w_{(2)}))
\]
for $v\in V$, $\lambda \in (W_1\otimes W_2)^*$,
$w_{(1)}\in W_1$ and $w_{(2)}\in W_2$.
For every  $\lambda \in (W_1\otimes W_2)^*$, the lower
truncation condition and the Jacobi identity clearly hold, and
$\tau_{P(z)}(Y_{t}({\bf 1},x))\lambda = \lambda$.  Hence the whole
(finite-dimensional) space
$(W_1\otimes W_2)^*$ is a $V$-module, and $W_1\warning W_2=(W_1\otimes
W_2)^*$.  This is in fact just the contragredient
module of the $V$-module $W_1\otimes W_2$ defined in Example
\ref{badim}.  We know from Example \ref{badim} that the identity map
on $(W_1\otimes W_2)^*$ is a quasi-$P(z)$-intertwining map that is not
in general a $P(z)$-intertwining map.  In Theorem \ref{pcorresp}, take
$F$ to be this identity map, so that $F^\vee$ is the identity map on
$(W_1\otimes W_2)^*$.  The proof of Theorem \ref{pcorresp} immediately
shows that $\lambda\in(W_1\otimes W_2)^*$ satisfies the
$P(z)$-compatibility condition if and only if
\begin{eqnarray*}
&{\displaystyle z^{-1}\delta\left(\frac{x_{1}-x_{0}}{z} \right)
\lambda((v \cdot w_{(1)}) \otimes w_{(2)})
+x^{-1}_{0}\delta\left(\frac{z-x_{1}}{-x_{0}}\right)
\lambda(w_{(1)} \otimes (v \cdot w_{(2)}))}&\nno\\
&{\displaystyle =\delta\left(\frac{x_{1}-z}{x_{0}} \right)
\lambda(w_{(1)} \otimes (v \cdot w_{(2)})),}&
\end{eqnarray*}
or equivalently,
\[
\lambda((v\cdot w_{(1)})\otimes w_{(2)})=\lambda(w_{(1)}\otimes(v\cdot
w_{(2)}))
\]
for any $v\in V$, $w_{(1)}\in W_1$ and $w_{(2)}\in W_2$, which of
course is not true in general (cf. Example \ref{badim}).  That is,
$W_1\warning W_2$ typically has a lot of non-$P(z)$-compatible
elements.  In fact, the space of compatible elements in $(W_1\otimes
W_2)^*$ naturally identifies with $(W_1\otimes_V W_2)^*$, the dual
space of the tensor product, over the commutative associative algebra
$V$, of the $V$-modules $W_1$ and $W_2$.  This space of compatible
elements is naturally a $V$-module (with $V$ viewed either as a
commutative associative algebra or as a vertex operator algebra)---
the contragredient module (with $V$ viewed either way) of
$W_1\otimes_V W_2$, which of course is naturally a quotient space of
$W_1\otimes W_2$.
In particular, take $W_1$ and $W_2$ to be the commutative
associative algebra $V$ itself, viewed as a module.
Then $V\warning V=(V\otimes V)^*$,
while the space of compatible elements $\lambda$ is naturally
identified with $(V \otimes_V V)^* = V^*$.
The contragredient module of $(V\otimes V)^*$ cannot equal the correct
tensor product module, namely, $V$, unless $V$ is $1$-dimensional.}
\end{rema}

\setcounter{equation}{0}
\setcounter{thm}{0}
\section{Further examples and $g(V)_{\geq 0}$-homomorphisms}

This section and the next are independent of Section 2, except for the
definition of the notion of (logarithmic) quasi-intertwining operator.

In this section we will give further examples of
quasi-intertwining operators that are not intertwining operators, by
using a canonical Lie algebra associated with a vertex
algebra and modules for this Lie algebra.  Just as in Section 2, such
examples give examples of non-$P(z)$-compatible modules.  Motivated by these
examples, we study some further properties of modules for this
canonical Lie algebra, which will lead us to other results about
quasi-intertwining operators.

Let $(V,Y,{\bf 1})$ be a vertex algebra.  Recall the canonical
Lie algebra $g(V)$ associated with $V$ (see \cite{B},
\cite{FFR}, \cite{Li}, \cite{MP}):
$$g(V)=(V\otimes \C[t,t^{-1}])/{\rm Im\;} ({\mathcal{D}}\otimes 1+1\otimes d/dt),$$
where ${\mathcal D}$ is given by ${\mathcal D}u=u_{-2}{\bf 1}
$ for
$u\in V$, with the bracket defined by means of representatives by:
\[
[u\otimes t^{m},v\otimes t^{n}]
=\sum_{i\ge 0}{m\choose i} (u_{i}v\otimes t^{m+n-i}).
\]

Denote by $\pi$ the natural quotient map
\[
\pi: V\otimes \C[t,t^{-1}]\rightarrow g(V).
\]
For $v\in V,\; n\in \Z$, we set
\[
v(n)=\pi (v\otimes t^{n})\in g(V).
\]
Note that ${\bf 1}(-1)$ is a nonzero central element of $g(V)$.
A $g(V)$-module
on which ${\bf 1}(-1)$ acts as a scalar $\ell$
is said to be of {\em level} $\ell$.
Recall that as in the case of affine Lie algebras, a
$g(V)$-module $W$ is said to be {\em restricted} if for any $w\in W$
and $v\in V$, we have $v(n)w=0$ for $n$ sufficiently large.
A $V$-module is automatically a restricted $g(V)$-module of level $1$
on which $v(n)$ acts as $v_n$. (See \cite{B}, \cite{FFR}, \cite{Li},
\cite{MP} for these and other standard notations and properties of
$g(V)$ and its modules.)

We set
\begin{eqnarray*}
& &g(V)_{\ge 0}=\pi(V\otimes \C[t])\subset g(V),\\
& &g(V)_{<0}=\pi(V\otimes t^{-1}\C[t^{-1}])\subset g(V).
\end{eqnarray*}
Clearly, these are Lie subalgebras of $g(V)$ and
\[
g(V)=g(V)_{<0}\oplus g(V)_{\ge 0}.
\]

The following observation, due to \cite{dlm-vpa} and \cite{Primc}, will be used
in the Appendix:

\bp{ldlm-p}
The linear map from $V$ to $g(V)_{<0}$ sending $v$ to
$v(-1)\;(=\pi(v\otimes t^{-1}))$
is a linear isomorphism. \epf
\ep

In addition to the examples in Example \ref{badim} (and Remark
\ref{rbadio}), we now give another
way of constructing examples of quasi-intertwining
operators that are not intertwining operators. We have:
\bp{pexample}
Let $V$ be a vertex operator algebra, let $(W_{1}, Y_{1})$ and
$(W_{2}, Y_{2})$ be weak $V$-modules, and
let $\theta$ be a
linear map from $W_{1}$ to $W_{2}$.
We define a linear map $\Y$ from $W_{1}$ to $\Hom (V,W_{2}((x)))$ by
\[
\Y (w,x)v=e^{xL(-1)}Y_{2}(v,-x)\theta(w)\;\;\;\mbox{ for }w\in W_{1},\; v\in V.
\]
Then $\Y$ is a quasi-intertwining operator (of type ${W_2 \choose W_1\,V}$)
if and only if $\theta$ is a
$g(V)_{\ge 0}$-homomorphism.
Furthermore, $\Y$ is an intertwining operator
if and only if $\theta$ is a $V$-homomorphism.
\ep
\pf Define
\[
Y_{2}^{t}(w,x)v=e^{xL(-1)}Y_{2}(v,-x)w\;\;\;\mbox{ for }w\in W_{2},\; v\in V.
\]
Then we have
\[
\Y(w,x)v=Y_{2}^{t}(\theta(w),x)v\;\;\;\mbox{ for }w\in W_{1},\; v\in V.
\]
{}From \cite{FHL}, $Y_{2}^{t}$ is an intertwining operator of type
${W_{2}\choose W_{2}\,V}$.
Assume that $\theta$ is a $g(V)_{\ge 0}$-homomorphism.
For any
$V$-module $(W, Y_{W})$ and any $u\in V$, we use
$Y_{W}^{-}(u,x)$ to denote $\sum_{n\ge 0}u_{n}x^{-n-1}$,
where $Y_{W}(u, x)=\sum_{n\in \Z}u_{n}x^{-n-1}$.
For any $u,v\in V, \; w\in W_{1}$,
using the fact that $\theta$ is a $g(V)_{\ge 0}$-homomorphism
we have
\begin{eqnarray}\label{equasi-half-hom}
\lefteqn{Y_{2}(u,x_{1})\Y(w,x_{2})v-\Y(w,x_{2})Y(u,x_{1})v}\nonumber\\
&&=Y_{2}(u,x_{1})Y_{2}^{t}(\theta(w),x_{2})v
-Y_{2}^{t}(\theta(w),x_{2})Y(u,x_{1})v\nonumber\\
&&=\Res_{x_{0}}x_{2}^{-1}\delta\left(\frac{x_{1}-x_{0}}{x_{2}}\right)
Y_{2}^{t}(Y_{2}(u,x_{0})\theta(w),x_{2})v\nonumber\\
&&=\Res_{x_{0}}x_{2}^{-1}\delta\left(\frac{x_{1}-x_{0}}{x_{2}}\right)
Y_{2}^{t}(Y_{2}^{-}(u,x_{0})\theta(w),x_{2})v\nonumber\\
&&=\Res_{x_{0}}x_{2}^{-1}\delta\left(\frac{x_{1}-x_{0}}{x_{2}}\right)
Y_{2}^{t}(\theta(Y_{1}^{-}(u,x_{0})w),x_{2})v\nonumber\\
&&=\Res_{x_{0}}x_{2}^{-1}\delta\left(\frac{x_{1}-x_{0}}{x_{2}}\right)
Y_{2}^{t}(\theta(Y_{1}(u,x_{0})w),x_{2})v\nonumber\\
&&=\Res_{x_{0}}x_{2}^{-1}\delta\left(\frac{x_{1}-x_{0}}{x_{2}}\right)
\Y(Y_{1}(u,x_{0})w,x_{2})v.
\end{eqnarray}
For $w\in W_{1}$, noticing that $\theta(L(-1)w)=L(-1)\theta(w)$ we have
\begin{eqnarray*}
& &\Y(L(-1)w,x)=Y_{2}^{t}(\theta(L(-1)w),x)=Y_{2}^{t}(L(-1)\theta(w),x)
=\frac{d}{dx}Y_{2}^{t}(\theta(w),x)\nonumber\\
& &=\frac{d}{dx}\Y(w,x).
\end{eqnarray*}
Thus $\Y$ is a quasi-intertwining operator.

Conversely, assume that $\Y$ is a quasi-intertwining operator.
Then the outside equality of (\ref{equasi-half-hom}) holds.
{}Using the first three  and the last two equalities in
(\ref{equasi-half-hom}) we see that
\begin{eqnarray*}
\lefteqn{\Res_{x_{0}}x_{2}^{-1}\delta\left(\frac{x_{1}-x_{0}}{x_{2}}\right)
Y_{2}^{t}(Y_{2}^{-}(u,x_{0})\theta(w),x_{2})v}\nonumber\\
&&=\Res_{x_{0}}x_{2}^{-1}\delta\left(\frac{x_{1}-x_{0}}{x_{2}}\right)
Y_{2}^{t}(\theta(Y_{1}^{-}(u,x_{0})w),x_{2})v.
\end{eqnarray*}
For any $n\ge 0$, we have
\begin{eqnarray*}
\lefteqn{\Res_{x_{0}}x_{2}^{-1}\delta\left(\frac{x_{1}-x_{0}}{x_{2}}\right)
x_{0}^{n}Y_{2}^{t}(Y_{2}^{-}(u,x_{0})\theta(w),x_{2})v}\nonumber\\
&&=\Res_{x_{0}}x_{2}^{-1}\delta\left(\frac{x_{1}-x_{0}}{x_{2}}\right)
(x_{1}-x_{2})^{n}Y_{2}^{t}(Y_{2}^{-}(u,x_{0})\theta(w),x_{2})v\nonumber\\
&&=\Res_{x_{0}}x_{2}^{-1}\delta\left(\frac{x_{1}-x_{0}}{x_{2}}\right)
(x_{1}-x_{2})^{n}Y_{2}^{t}(\theta(Y_{1}^{-}(u,x_{0})w),x_{2})v\nonumber\\
&&=\Res_{x_{0}}x_{2}^{-1}\delta\left(\frac{x_{1}-x_{0}}{x_{2}}\right)
x_{0}^{n}Y_{2}^{t}(\theta(Y_{1}^{-}(u,x_{0})w),x_{2})v.
\end{eqnarray*}
Taking $\Res_{x_{1}}$ we get
\[
Y_{2}^{t}(u_{n}\theta(w), x_{2})v=Y_{2}^{t}(\theta(u_{n}w),x_{2})v.
\]
Setting $v={\bf 1}$ in this formula and using the definition of $Y_{2}^{t}$,
we obtain $u_{n}\theta(w)=\theta(u_{n}w)$.
This proves that $\theta$ is a $g(V)_{\ge 0}$-homomorphism.

Using the whole Jacobi identity instead of
the commutator formula one shows analogously that $\Y$ is an intertwining operator
if and only if $\theta$ is a $g(V)$-homomorphism.
\epfv

Partly due to this proposition, we are interested in looking for
$g(V)_{\geq 0}$-module maps that are not $V$-module maps.  We will study this
problem in the general context of vertex algebras.  We now fix a
vertex algebra $(V,Y,{\bf 1})$.

Recall the following result (\cite{DLM1}, \cite{LL}, Proposition
4.5.7):

\bp{p-ll}
Let $W$ be a $V$-module and
let $u,v\in V,\; p,q\in \Z$ and $w\in W$.
Let $l$ be a nonnegative integer such that
$$u_{n}w=0\;\;\;\mbox{ for }n\ge l$$
and let $m$ be a nonnegative integer such that
$$v_{n}w=0\;\;\;\mbox{ for }n>m+q.$$
Then
\begin{eqnarray}\label{eassociativity-formula}
u_{p}v_{q}w=\sum_{i=0}^{m}\sum_{j=0}^{l}{p-l\choose i}{l\choose j}
(u_{p-l-i+j}v)_{q+l+i-j}w.
\end{eqnarray}
\epf
\ep

The following result, due to \cite{DM} and \cite{Li}, is an
immediate consequence of Proposition \ref{p-ll}:

\bp{ldm-li}
Let $W$ be a $V$-module.
Then for any $w\in W$,
the linear span of the vectors
$v_{n}w$ for $v\in V$, $n\in \Z$,
namely, $g(V)w$, is a
$V$-submodule of $W$. Furthermore,
for any subspace or subset $U$ of $W$, $g(V)U$ is the
$V$-submodule of $W$
generated by $U$. \epf
\ep

Examining Proposition \ref{p-ll} more closely we have:

\bp{gvsub}
Let $W$ be a $V$-module.
Then $g(V)_{\ge 0}w$ is a
$V$-submodule of $W$ for every $w\in W$, and
so is $g(V)_{\ge 0}W$.
Furthermore, $g(V)_{\ge 0}V$ is a (two-sided) ideal of $V$.
\ep
\pf It follows directly from the formula
(\ref{eassociativity-formula}) with $q\ge 0$ that for any $w\in W$,
$g(V)_{\ge 0}w$ is a
$V$-submodule of $W$.
Thus $g(V)_{\ge 0}W$ is a
$V$-submodule.
Taking $W=V$, we see that
$g(V)_{\ge 0}V$ is a left ideal of $V$.
Since $[{\mathcal{D}}, v_{n}]=-nv_{n-1}$ for $v\in V,\; n\in \Z$,
we have that $[{\mathcal{D}}, g(V)_{\ge 0}]\subset g(V)_{\ge 0}$,
acting on $V$.
Therefore ${\mathcal{D}} g(V)_{\ge 0}V\subset g(V)_{\ge 0}V$.
{}From Remark 3.9.8 in \cite{LL}, $g(V)_{\ge 0}V$ is an ideal.
\epfv

We have:

\begin{propo}\label{p-equivalence}
The following statements are equivalent:

\begin{enumerate}

\item There exists a
$V$-module $W$ such that
$g(V)_{\ge 0}W\ne W$.

\item  $g(V)_{\ge 0}V\ne V$, or
equivalently ${\bf 1}\notin g(V)_{\ge 0}V$.

\end{enumerate}
\end{propo}

\pf  Since $g(V)_{\ge 0}V$ is an ideal (Proposition \ref{gvsub}), the two
conditions in the second statement are equivalent.
We need only prove that if there exists a
$V$-module $W$
such that $g(V)_{\ge 0}W\ne W$, then $g(V)_{\ge 0}V\ne V$.
Let $W$ be such a
$V$-module.
Then we have a nonzero
$V$-module $\widetilde{W}=W/g(V)_{\ge 0}W$ (see Proposition \ref{gvsub})
such that $g(V)_{\ge 0}\widetilde{W}=0$.
Since $Y({\bf 1},x)=1$ on $\widetilde{W}$, the annihilating ideal
of $V$
\[
{\rm Ann}_{V}(\widetilde{W})=\{ v\in V\;|\; Y(v,x)\widetilde{W}=0\}
\]
(see Proposition 4.5.11 of \cite{LL}) is proper. For $u,v\in V,\; n\ge 0$, we have
\[
Y(u_{n}v,x)=\Res_{x_{1}}(x_{1}-x)^{n}[Y(u,x_{1}),Y(v,x)]=\sum_{i=0}^{n}
{n\choose i}(-x)^{i}[u_{n-i},Y(v,x)]
\]
(recall (3.8.14) in \cite{LL}).
Since $g(V)_{\ge 0}\widetilde{W}=0$, it follows that $u_{n}v\in {\rm
Ann}_{V}(\widetilde{W})$ for $u,v\in V,\; n\ge 0$. This proves that
$g(V)_{\ge 0}V\subset {\rm Ann}_{V}(\widetilde{W})$, a proper
subspace of $V$. Consequently, $g(V)_{\ge 0}V\ne V$.
\epfv

Furthermore, we have:

\bp{lnonzero-cc}
Suppose that $(V,Y,{\bf 1},\omega)$ is a vertex operator algebra of
central charge $c\ne 0$.
Then $g(V)_{\ge 0}W=W$ for any weak $V$-module $W$.
In particular, $g(V)_{\ge 0}V=V$.
\ep
\pf We have the following relation in $V$:
\[
L(2)\omega=L(2)L(-2){\bf 1}=\frac{1}{2}c{\bf 1}.
\]
Since $c\ne 0$, we have ${\bf 1}=(2/c)L(2)\omega\in g(V)_{\ge 0}V$.
By Proposition \ref{p-equivalence},
$g(V)_{\ge 0}W=W$ for any weak $V$-module $W$. \epfv

The following example shows that in Proposition \ref{lnonzero-cc}, the
condition $c\ne 0$ is necessary:
\begin{exam}\label{exa2.4}{\rm
Let $V$ be the minimal vertex operator algebra $V_{Vir}(0,0)$ associated with
the Virasoro algebra ${\mathcal{L}}$ of central charge $c=0$ (cf.
Section 6.1 of \cite{LL}).
We are going to show that $V_{(0)}\not\subset g(V)_{\ge 0}V$.
We know that $V_{(0)}=\C {\bf 1}$, $V_{(1)}=0$ and $V_{(2)}=\C \omega$,
where $\omega$ is the conformal vector.
Set $V_{+}=\coprod_{n\ge 1}V_{(n)}$.  We will show that $V_{+}$ is an ideal.
If this is proved, we will have
$$g(V)_{\ge 0}V=g(V)_{\ge 0}{\bf 1}+g(V)_{\ge 0}V_{+}=g(V)_{\ge
0}V_{+}\subset V_{+},$$
which immediately implies that ${\bf 1}\notin g(V)_{\ge 0}V$.
Note that from Section 6.1 of \cite{LL}, $U({\mathcal{L}})\omega$ is
a left ideal of $V$ (and hence, by Remark 3.9.8 of \cite{LL}, a
(two-sided) ideal of $V$).
It suffices to prove that $V_{+}=U({\mathcal{L}})\omega$.
Since $V$ is an ${\mathcal{L}}$-module with ${\bf 1}$ as a generator,
$V_{+}$ is spanned by the vectors
$$L(-m_{1})\cdots L(-m_{k}){\bf 1}$$
for $k\ge 1,\; m_{1}\ge \cdots \ge m_{k}\ge 2$.
Using the formula $Y(\omega,x){\bf 1}=e^{xL(-1)}\omega$
we have
$$L(n){\bf 1}\in U({\mathcal{L}})\omega\;\;\;\mbox{ for }n\in \Z.$$
It follows that $V_{+}\subset U({\mathcal{L}})\omega$.
On the other hand, for $n\ge 1$ we have
$$L(n)\omega=L(n)L(-2){\bf 1}=(n+2)L(n-2){\bf
1}+\delta_{n,2}\frac{1}{2}c{\bf 1}=0,$$
since the central charge $c$ is zero. That is, $\omega$ is a lowest
weight vector for the Virasoro algebra ${\mathcal{L}}$.  Thus
$$U({\mathcal{L}})\omega \subset V_{+}.$$
Therefore we have $V_{+}=U({\mathcal{L}})\omega$, proving that $V_{+}$
is an ideal of $V$, and hence proving that $V_{(0)}\not\subset
g(V)_{\ge 0}V$.  Using an analogous argument we easily also see that for
$V=V_{\hat{\g}}(0,0)$, associated with an affine Kac-Moody Lie algebra
$\hat{\g}$ of level $0$ (cf. Section 6.2 of \cite{LL}),
we have ${\bf 1}\notin g(V)_{\ge 0}V$. }
\end{exam}

We have:
\bp{phomomorphism}
Let $V$ be a vertex algebra and let $W_1$ and $W_2$ be
$V$-modules.
If $g(V)_{\ge 0}W_{1}=W_{1}$,
then any $g(V)_{\ge 0}$-homomorphism from
$W_{1}$ to $W_{2}$ is a $V$-homomorphism.
\ep
\pf Let $f$ be a $g(V)_{\ge 0}$-homomorphism from
$W_{1}$ to $W_{2}$. Let $u,v\in V,\; p,q\in \Z,\; w\in W_1$.
Let $l,m$ be nonnegative integers such that
\begin{eqnarray*}
& &u_{n}w=0,\;\;\; u_{n}f(w)=0\;\;\;\mbox{ for }n\ge l,\\
& &v_{n}w=0,\;\;\; v_{n}f(w)=0\;\;\;\mbox{ for }n>m+q.
\end{eqnarray*}
By Proposition \ref{p-ll} we have
\begin{eqnarray*}
& &u_{p}v_{q}w=\sum_{i=0}^{m}\sum_{j=0}^{l}{p-l\choose i}{l\choose j}
(u_{p-l-i+j}v)_{q+l+i-j}w,\label{eupvqw1}\\
& &u_{p}v_{q}f(w)=\sum_{i=0}^{m}\sum_{j=0}^{l}{p-l\choose i}{l\choose j}
(u_{p-l-i+j}v)_{q+l+i-j}f(w).\label{eupvqfw2}
\end{eqnarray*}
Notice that if $q\ge 0$, we have $q+l+i-j\ge 0$ for $i\ge 0,\; 0\le j\le l$.
Then for $q\ge 0$, using (\ref{eupvqw1}) and (\ref{eupvqfw2}) and
the fact that $f$ is a $g(V)_{\ge 0}$-homomorphism, we get
\begin{eqnarray*}
f(u_{p}v_{q}w)&=&\sum_{i=0}^{m}\sum_{j=0}^{l}{p-l\choose i}{l\choose j}
f\left((u_{p-l-i+j}v)_{q+l+i-j}w\right)\nonumber\\
&=&\sum_{i=0}^{m}\sum_{j=0}^{l}{p-l\choose i}{l\choose j}
(u_{p-l-i+j}v)_{q+l+i-j}f(w)\nonumber\\
&=&u_{p}v_{q}f(w)\nonumber\\
&=&u_{p}f(v_{q}w).
\end{eqnarray*}
This shows that
$$f(u_{p}w')=u_{p}f(w')\;\;\;\mbox{ for }u\in V,\; p\in \Z,\; w'\in g(V)_{\ge 0}W_1.$$
Since we are assuming that $W_{1}=g(V)_{\ge 0}W_{1}$,
$f$ is a $V$-homomorphism.
\epf

\bex{exaple1}
{\rm Here we give concrete examples of $g(V)_{\ge 0}$-homomorphisms
between $V$-modules that are not $V$-homomorphisms, for suitable
vertex operator algebras $V$.  Then
by Proposition \ref{pexample} we obtain quasi-intertwining operators
that are not intertwining operators, and as in Section 2, this gives
examples of non-$P(z)$-compatible modules.
Let $V$ be the vertex operator algebra $V_{Vir}(0,0)$
associated to the Virasoro algebra of central charge $0$,
or the vertex operator algebra  $V_{\hat{\g}}(0,0)$
associated with an affine Lie algebra $\hat{\g}$ of level
$0$.  Recall from Example \ref{exa2.4} that
$V=V_{+}\oplus \C {\bf 1}$, where $V_{+}=\coprod_{n>0}V_{(n)}$
is an ideal of $V$.
Consequently, $V/V_{+}\simeq \C$ is a (one-dimensional)
nontrivial $V$-module, on which
$Y(v,x)$ acts as zero for $v\in V_{+}$ (and $Y({\bf 1},x)$ acts as the
identity).
Define a linear map
\[
\theta: \C \rightarrow V;\;\; \alpha \mapsto \alpha {\bf
1}\;\;\;\mbox{ for }\alpha \in \C.
\]
Clearly, $\theta$ is a $g(V)_{\ge 0}$-homomorphism.
But $\theta$ is not a $V$-homomorphism, since for any nonzero $v\in
V_+$, $v_{-1}$ acts on $\C$ as zero but $v_{-1}{\bf 1}=v\ne 0$.  }
\eex

\setcounter{equation}{0}
\setcounter{thm}{0}
\section{Necessary conditions for the existence of examples}

In this section we let $V$ be a vertex operator algebra. We will show
that for weak $V$-modules $W_1$, $W_2$ and $W_3$,
the condition $W_1=g(V)_{\ge 0}W_1$
implies that any quasi-intertwining operator of type ${W_3\choose
W_1\,W_2}$ is an intertwining operator; in other words, the condition
$W_1 \ne g(V)_{\ge 0}W_1$ is necessary for the existence of a
quasi-intertwining operator of type ${W_3\choose W_1\,W_2}$ that is
not an intertwining operator.  We conclude that if $V$ has nonzero
central charge, then any quasi-intertwining operator among any weak
$V$-modules is an intertwining operator.

Recall the following definition from \cite{Li-local} and \cite{Li-hom}:

\bd{dhom}
{\rm Let
$(W_{2}, Y_{2}), (W_{3}, Y_{3})$ be
weak $V$-modules.
Denote by $\H(W_{2},W_{3})$
the vector subspace of $(\Hom (W_{2},W_{3}))\{x\}$ consisting
of the formal series $\phi(x)$ satisfying the following conditions:
Writing
\[
\phi(x)w_{(2)}=\sum_{n\in\C}w_{(3)}^{(n)}x^{-n-1}\;\;\;
(\mbox{where}\;\; w_{(3)}^{(n)} \in W_3)
\]
for $w_{(2)} \in W_2$, we have
\[
w_{(3)}^{(n)} = 0 \;\;  \mbox{for}\;\; n \;\; \mbox{whose real part is
sufficiently large;}
\]
\[
[L(-1), \phi(x)]=\frac{d}{dx}\phi(x);
\]
and for any $v\in V$, there exists a nonnegative
integer $k$ such that
\[
(x_{1}-x_{2})^{k}(Y_{3}(v,x_{1})\phi(x_{2})-\phi(x_{2})
Y_{2}(v,x_{1}))=0.
\]
We also define a vertex operator map
\[
Y_{\H}(\cdot,x_0):V \to (\End \H(W_{2},W_{3}))[[x_0,x_0^{-1}]]
\]
by
\begin{eqnarray}\label{YH}
\lefteqn{Y_{\H}(v,x_{0})\phi(x)}\nonumber\\
&&=\Res_{x_{1}}\left
( x_{0}^{-1}\delta\left(\frac{x_{1}-x}{x_{0}}\right)
Y_{3}(v,x_{1})\phi(x)-x_{0}^{-1}\delta\left(\frac{x-x_{1}}{-x_{0}}\right)
\phi(x)Y_{2}(v,x_{1})\right).
\end{eqnarray}}
\ed

The following was essentially proved in \cite{Li}:

\bt{told}
Let $W_{2}$ and $W_{3}$ be
weak $V$-modules.
Then $(\H(W_{2},W_{3}),Y_{\H})$ carries the structure of
a weak $V$-module. Furthermore,
for any weak $V$-module $W_{1}$,
giving an intertwining operator $\Y$ of type
${W_{3}\choose W_{1}W_{2}}$ is equivalent to giving
a $V$-homomorphism $\psi = \psi_{x}$
{}from $W_{1}$ to $\H(W_{2},W_{3})$, where
the relation between $\Y$ and $\psi$ is given by
\[
\psi_{x}(w_{(1)})=\Y(w_{(1)},x)\;\;\;\mbox{ for }w_{(1)}\in W_{1}.
\]
\epf
\et

\br{r-old}
{\rm For $V$-modules $W_{2},W_{3}$, the space $\H(W_{2},W_{3})$
was defined in \cite{Li-local}
and \cite{Li-hom}, where it was denoted by $G(W_{2},W_{3})$.
Theorem \ref{told} was proved in \cite{Li-local}
and \cite{Li-hom} with $W_{2},W_{3}$ being $V$-modules, but the proof did
not use the grading.}
\er

Let $W_{1}, W_{2}$ and $W_{3}$ be
$V$-modules
and  $\Y$ a quasi-intertwining operator of type ${W_{3}\choose
W_{1}W_{2}}$. Since $L(-1)=\omega_{0}$, by using the commutator formula
and the $L(-1)$-derivative property we get
\[
[L(-1),\Y(w_{(1)},x)]=\Y(L(-1)w_{(1)},x)
=\frac{d}{dx}\Y(w_{(1)},x)\;\;\;\mbox{ for }w_{(1)}\in W_{1}.
\]
For $v\in V,\; w_{(1)}\in W_{1}$, let $k$ be a nonnegative integer such that
$v_{n}w_{(1)}=0$ for $n\ge k$. Using the commutator formula we get
\[
(x_{1}-x_{2})^{k}[Y(v,x_{1}),\Y(w_{(1)},x_{2})]=0.
\]
Thus we have
\[
\Y(w_{(1)},x)\in \H (W_{2},W_{3})\;\;\;\mbox{ for }w_{(1)}\in W_{1}.
\]
In view of this, we may and should consider $\Y = \Y(\cdot,x)$ as a linear map from
$W$ to $\H(W_{2},W_{3})$.
Furthermore, for $v\in V,\; n\ge 0,\; w_{(1)}\in W_{1}$,
using the definition of the action of $v_{n}$ on $\H(W_{2},W_{3})$
and the commutator formula for $\Y$ we get
\begin{eqnarray*}
v_{n}(\Y(w_{(1)},x))
&=&\Res_{x_{1}}(x_{1}-x)^{n}[Y(v,x_{1}),\Y(w_{(1)},x)]\nonumber\\
&=&\Res_{x_{1}}\Res_{x_{0}}(x_{1}-x)^{n}
x_{1}^{-1}\delta\left(\frac{x+x_{0}}{x_{1}}\right)\Y(Y(v,x_{0})w_{(1)},x)\nonumber\\
&=&\Res_{x_{0}}x_{0}^{n}\Y(Y(v,x_{0})w_{(1)},x)\nonumber\\
&=&\Y(v_{n}w_{(1)},x).
\end{eqnarray*}
Thus $\Y$ is a $g(V)_{\ge 0}$-homomorphism. Therefore we have proved:

\bp{pquasi-IO}
Let $W_{1}, W_{2}$ and $W_{3}$ be weak $V$-modules
and $\Y$ be a quasi-intertwining operator of type
${W_{3}\choose W_{1}W_{2}}$. Then
\[
\Y(w_{(1)},x)\in \H(W_{2},W_{3})\;\;\;\mbox{ for }w_{(1)}\in W_{1}.
\]
Furthermore, the linear map $\psi_x$ from $W_{1}$ to $\H(W_{2},W_{3})$
defined by
\[
\psi_{x}(w_{(1)})=\Y(w_{(1)},x)
\]
is a $g(V)_{\ge 0}$-homomorphism, i.e.,
\[
\psi_{x} (v_{n}w_{(1)})=v_{n}\psi_{x}(w_{(1)})
\]
for $v\in V,\; n\ge 0,\; w_{(1)}\in W_{1}$. \epf
\ep

In view of Theorem \ref{told} and Proposition \ref{pquasi-IO},
Proposition \ref{phomomorphism} gives:

\bp{pquasi-int-new}
Let  $W_{1},W_{2}$ and $W_{3}$
be weak $V$-modules. If $g(V)_{\ge 0}W_1=W_1$,
then any quasi-intertwining
operator of type ${W_{3}\choose W_{1}W_{2}}$ is an
intertwining operator.
\ep

\pf Let $\Y$ be a quasi-intertwining operator of
type ${W_{3}\choose W_{1} W_{2}}$. By Proposition \ref{pquasi-IO},
there exists a $g(V)_{\ge 0}$-homomorphism
$\psi_x$ from $W_{1}$ to the weak $V$-module $\H (W_{2},W_{3})$ such that
\[
\Y(w,x)=\psi_{x} (w)\;\;\;\mbox{ for }w\in W_{1}.
\]
{}From Theorem \ref{told}, $\Y$ is an intertwining operator if and only if
$\psi_x$ is a $V$-homomorphism, i.e., a $g(V)$-homomorphism.
But by Proposition \ref{phomomorphism},
$\psi_x$ is indeed a $V$-homomorphism. Thus $\Y$ is an intertwining
operator. \epfv

Combining Propositions \ref{p-equivalence},
\ref{phomomorphism} and \ref{pquasi-int-new},
we immediately have the first assertion of the following theorem:

\bt{thalf-homomorphism-general}
Suppose that $g(V)_{\ge 0}V= V$.
Then any $g(V)_{\ge 0}$-homomorphism between weak $V$-modules is a
$V$-homomorphism and any quasi-intertwining
operator among any weak $V$-modules is an
intertwining operator.
On the other hand, if $g(V)_{\ge 0}V\ne V$ and $\dim V>1$, then there exists
a $g(V)_{\ge 0}$-homomorphism
between $V$-modules that is not a $V$-homomorphism
and there exists a quasi-intertwining
operator among $V$-modules that is not an
intertwining operator.
\et

\pf Assume $g(V)_{\ge 0}V\ne V$ and $\dim V>1$. We modify the construction
given in Example \ref{exaple1} as follows.
Set $W=V/g(V)_{\ge 0}V$. Since $g(V)_{\ge 0}V$ is an ideal of $V$
(by Proposition \ref{gvsub}), $W$ is a nonzero
module for the vertex operator algebra $V$, and we have
$g(V)_{\ge 0}W=0$. Let $f$ be any nonzero linear functional on $W$.
Define
$$\theta_{f}: W\rightarrow V,\;\; w\mapsto f(w){\bf 1}.$$
With $g(V)_{\ge 0}W=0$ and $g(V)_{\ge 0}{\bf 1}=0$,
it is clear that $\theta_{f}$ is a $g(V)_{\ge 0}$-homomorphism.
Let $v\in V\backslash\C {\bf 1}$ (as $\dim V>1$) and let $w\in W$ be such that $f(w)=1$.
We have \begin{eqnarray*}
& &v_{-1}\theta_{f}(w)=v_{-1}f(w){\bf 1}=v_{-1}{\bf 1}=v\notin \C{\bf 1},\\
& &\theta_{f}(v_{-1}w)=f(v_{-1}w){\bf 1}\in \C{\bf 1},
\end{eqnarray*}
proving that $\theta_{f}$ is not a $V$-homomorphism. By Proposition
\ref{pexample} this yields to a quasi-intertwining operator (of
type ${V\choose W V}$) that is not an
intertwining operator.
\epf

\br{r1-dimension}
{\em If $\dim V=1$, i.e., $V=\C {\bf 1}$ with ${\bf 1}\ne 0$,
then $g(V)_{\ge 0}V=0\ne V$. In this case, a $V$-homomorphism
is simply a linear map, so that
any $g(V)_{\ge 0}$-homomorphism is a $V$-homomorphism.}
\er

Combining Proposition \ref{lnonzero-cc} with
Theorem \ref{thalf-homomorphism-general}
we immediately have:

\bc{chalf-homomorphism}
Suppose that the central charge of $V$ is not $0$.
Then any $g(V)_{\ge 0}$-homomorphism between weak $V$-modules is a
$V$-homomorphism and any quasi-intertwining
operator among any weak $V$-modules is an
intertwining operator. \epf
\ec

\begin{rema}\label{Hlog}
{\rm
Let $W_{2}$ and $W_{3}$ be weak $V$-modules.
By analogy with the space $\H(W_2, W_3)$, define
$\H^{\log}(W_2,W_3)$ to be the vector subspace of $(\Hom
(W_{2},W_{3}))\{x\}[[\log x]]$ consisting of the formal series $\phi(x)$ satisfying
the following conditions: Writing
\[
\phi(x)w_{(2)}=\sum_{n\in\C}\sum_{l\in\N}w_{(3)}^{(n;\,l)}
x^{-n-1}(\log x)^{l}\;\;\;
(\mbox{where}\;\; w_{(3)}^{(n;\,l)} \in W_3)
\]
for $w_{(2)} \in W_2$, we have
$$w_{(3)}^{(n;\,l)} = 0$$
if either $l\in {\Bbb N}$ is sufficiently large or the real part of
$n$ is sufficiently large;
$$[L(-1), \phi(x)]=\frac{d}{dx}\phi(x);$$
and for any $v\in V$, there exists a nonnegative
integer $k$ such that
$$(x_{1}-x_{2})^{k}[Y(v,x_{1}), \phi(x_{2})]=0.$$
Define a vertex operator map $Y_{\H}$ from $V$ to
$(\End \H^{\log}(W_{2},W_{3}))[[x, x^{-1}]]$ by (\ref{YH}).
Note that the coefficients of the
powers of $\log x$ in logarithmic intertwining operators (recall
Definition \ref{logio}) satisfy the Jacobi identity.
Following the proof of
Theorem \ref{told} in \cite{Li}, we have that the coefficients of the
powers of $\log x$ in elements of $\H^{\log}(W_2,W_3)$
satisfy the Jacobi identity and thus
$(\H^{\log}(W_{2},W_{3}),Y_{\H})$ carries
the structure of a weak $V$-module.
Then the same proof as that for Proposition
\ref{pquasi-IO} shows that a logarithmic
intertwining operator $\Y$ of type ${W_{3}\choose W_{1}W_{2}}$
gives a natural $V$-homomorphism $\psi$ from $W_{1}$ to
$\H^{\log}(W_{2},W_{3})$ and that a logarithmic quasi-intertwining
operator of type ${W_{3}\choose W_{1}W_{2}}$ (recall Definition
\ref{logquasi}) gives a natural
$g(V)_{\ge 0}$-homomorphism from $W_{1}$ to $\H^{\log}(W_{2},W_{3})$.
Hence we see that the statements of
Proposition \ref{pquasi-int-new}, Theorem
\ref{thalf-homomorphism-general} and Corollary \ref{chalf-homomorphism}
also hold with ``quasi-intertwining operator''
replaced by ``logarithmic quasi-intertwining operator,'' and
``intertwining operator'' replaced by ``logarithmic intertwining
operator.''  }
\end{rema}

\renewcommand{\theequation}{A.\arabic{equation}}
\renewcommand{\thethm}{A.\arabic{thm}}
\setcounter{equation}{0}
\setcounter{thm}{0}
\section*{Appendix: The Jacobi identity vs.~the commutator formula
for modules}

In the definition of the notion of module for a vertex (operator)
algebra, is the commutator formula enough?  That is, does the
commutator formula (see (\ref{commutatorformula}) below)
imply the Jacobi identity? The answer is no, as one
would expect. In fact the easiest counterexample is quite simple.  The
following is taken from Remark 4.4.6 of \cite{LL}:

\begin{exam}{\rm
Let $V$ be the $2$-dimensional commutative associative algebra with
a basis $\{1,a\}$ such that $a^{2}=1$.
Then $V$ has a vertex operator algebra structure with $Y(u,x)v=uv$
for $u,v\in V$ and with ${\bf 1}=1$ and $\omega = 0$.
Equip the $1$-dimensional space $W=\C w$ with a
linear map $Y_W: V\to\Hom(W,W((x)))$ determined by $Y_W(1,x)w=w$,
$Y_W(a,x)w=0$.  Then $(W,Y_W)$ satisfies all
the axioms for a $V$-module except the Jacobi identity, and the
commutator formula certainly holds (trivially).  In fact the
Jacobi identity fails since $Y_W(Y(a,x_0)a,x_2)w=w \neq
0=Y_W(a,x_0+x_2)Y_W(a,x_2)w$.  }
\end{exam}

We now give some less trivial counterexamples.

Let $V$ be a vertex operator algebra.
Let $(W,Y_W)$ be a pair that satisfies all
the axioms in the definition of the notion of module for $V$ viewed as
a vertex
algebra (see Definition 4.1.1 in \cite{LL}) except that the Jacobi
identity is replaced by the commutator formula:
\begin{eqnarray}\label{commutatorformula}
[Y_{W}(u,x_{1}),Y_{W}(v,x_{2})]
=\Res_{x_{0}}x_{2}^{-1}\delta\left(\frac{x_{1}-x_{0}}{x_{2}}\right)
Y_{W}(Y(u,x_{0})v,x_{2})
\end{eqnarray}
for $u,v\in V$, and in addition assume that the $L(-1)$-derivative
property also holds on $W$. Then $W$ is naturally a restricted
$g(V)$-module of level $1$.  Conversely, let $W$ be a restricted
$g(V)$-module of level $1$. Define a linear map $Y_{W}: V\rightarrow
\Hom (W,W((x)))$ by
\[
Y_W(v,x)=\sum_{n\in\Z}v(n)x^{-n-1}.
\]
Then $(W, Y_W)$ satisfies all the axioms in the definition of weak
$V$-module except that the Jacobi identity is replaced by the
commutator formula, and in addition, the $L(-1)$-derivative property holds.
We have:
\begin{propo}\label{badmodule}
Unless the vertex operator algebra
$V$ is one-dimensional, there exists a restricted $g(V)$-module
of level $1$ that is not a weak $V$-module. Furthermore, such an
example can be chosen to satisfy the two grading restriction
conditions in the definition of the notion of $V$-module
if $V$ has no elements of negative weight and $V_{(0)}=\C {\bf 1}$.
In particular, for any such vertex operator algebra $V$, the Jacobi
identity cannot be replaced by the commutator formula in the
definition of the notion of module.
\end{propo}

\pf In view of the creation property and vacuum property,
$\C{\bf 1}$ is a $(g(V)_{\ge 0}\oplus \C{\bf 1}(-1))$-submodule of
$V$, with $g(V)_{\ge 0}$ acting trivially and ${\bf 1}(-1)$ acting
as the identity.
Form the induced $g(V)$-module
$$W=U(g(V))\otimes _{U(g(V)_{\ge 0}\oplus \C {\bf 1}(-1))}\C {\bf 1},$$
which is of level $1$.
It follows by induction that $W$ is a restricted $g(V)$-module (of level $1$).
By the Poincar\'e-Birkhoff-Witt theorem and Proposition \ref{ldlm-p} we have
\begin{eqnarray}\label{elast}
W=U(g(V)_{<0}/\C{\bf 1}(-1)) \simeq S(V/\C {\bf 1})
\end{eqnarray}
as a vector space. Notice that $\wt v_{-n}=\wt v+n-1>0$ for homogeneous vector $v$
of positive weight and for $n\ge 1$. If $V$
has no elements of negative weight and $V_{(0)}=\C {\bf 1}$, $W$
satisfies the two grading restriction conditions
in the definition of the notion of $V$-module.

Now, we claim that $W$ is not a weak $V$-module if $\dim V>1$.
Otherwise, with $g(V)_{\ge 0}(1\otimes {\bf 1})=0$, the standard
generator $1\otimes {\bf 1}$ of $W$ is a vacuum-like vector and we
have a $V$-homomorphism from $V$ into $W$ sending $v$ to
$v_{-1}(1\otimes {\bf 1})\;(=v(-1)\otimes {\bf 1})$ for $v \in V$ (see
\cite{Li-form}; cf. Section 4.7 of \cite{LL}).
The image of this map is $g(V)_{< 0}(1\otimes {\bf 1}) = g(V)(1\otimes
{\bf 1})$ by Proposition \ref{ldlm-p}, and this space is a
$V$-submodule of $W$ by Proposition \ref{ldm-li}.  Thus the map is surjective.
That is, $$W=\{ v(-1)\otimes {\bf 1} \;|\; v\in V\}.$$
But (\ref{elast}) implies that $\{ v(-1)\otimes {\bf 1} \;|\; v\in V\}$
is a proper subspace of $W$ unless $V/\C {\bf 1}=0$.
Thus $W$ is not a weak $V$-module if $\dim V>1$.
\epfv

For a vertex operator algebra $V$, the Lie algebra $g(V)$
is naturally a $\Z$-graded Lie algebra $g(V)=\coprod_{n\in \Z}g(V)_{(n)}$,
where the $\Z$-grading is given by $L(0)$-weights:
\[
\wt (u\otimes t^{m})=\wt u-m-1
\]
for homogeneous $u$ and for $m\in \Z$.
For any $n\in \Z$, we have
\begin{eqnarray*}
g(V)_{(n)}&=&\sum_{m\in \Z}\pi(V_{(m)}\otimes t^{m-1-n})\\
&=&{\rm span}\{ v(m-1-n)\;|\; v\in V_{(m)},\; m\in \Z\}.
\end{eqnarray*}
We set
\[
g(V)_{(\pm)}=\coprod_{n>0}g(V)_{(\pm n)}.
\]
(Note the distinction between $g(V)_{(-)}$ and $g(V)_{<0}$.)
Clearly, the Lie subalgebras $g(V)_{\ge 0}$ and $g(V)_{<0}$ are also
graded subalgebras of $g(V)$.
Note that if $V$ has no elements of negative weight, i.e.,
if $V=\coprod_{n\ge 0}V_{(n)}$, then
$g(V)_{(-)}$ is a subalgebra of $g(V)_{\ge 0}$.
If in addition $V_{(0)}=\C {\bf 1}$, we also have
$g(V)_{(0)}\subset g(V)_{\ge 0}\oplus \C {\bf 1}(-1)$.

\begin{rema}{\rm
Here we give another construction of counterexamples.
Let $V$ be any nonzero vertex operator algebra.
If the conformal vector $\omega$ is zero, then
$V={\rm Ker \,} L(-1)=V_{(0)}$ is simply a finite-dimensional
commutative associative algebra with identity.
If $\dim V=1$, a weak $V$-module simply amounts to a vector space.
If $\dim V>1$, we have already seen that
a restricted $g(V)$-module of level $1$ is not necessarily a $V$-module.
Now assume that $\omega\ne 0$. Let $V=\coprod_{n\ge 0}V_{(r+n)}$
with $V_{(r)}\ne 0$.
Then $V_{(r)}$ is naturally a module for the Lie subalgebra
 $g(V)_{(0)}\oplus g(V)_{(-)}$ of $g(V)$.
Consider the induced $g(V)$-module
$$M=U(g(V))\otimes _{U(g(V)_{(0)}\oplus g(V)_{(-)})}V_{(r)}.$$
Clearly, $M$ is a $\Z$-graded $g(V)$-module of level $1$ where
the grading is given by the $L(0)$-eigenspaces.
By the Poincar\'e-Birkhoff-Witt theorem,
\begin{eqnarray}\label{eM-PBW}
M \simeq U(g(V)_{(+)})\otimes V_{(r)},
\end{eqnarray}
which implies that $M_{(n)}=0$ for $n<r$ and $M_{(r)}=V_{(r)}$.
Consequently, $M$ is a restricted $g(V)$-module (of level $1$) and
hence $M$ satisfies all the axioms in the definition of weak
$V$-module except that the Jacobi identity is replaced by the
commutator formula.
We claim that $M$ is not a weak $V$-module.
Otherwise, by Proposition \ref{ldm-li} we would have
$$M=g(V)V_{(r)}=g(V)_{(+)}V_{(r)}+V_{(r)}.$$
Combining this with (\ref{eM-PBW}) we must have $g(V)_{(+)}=0$.
But $0\ne L(-2)\in g(V)_{(+)}$, since
$L(-2){\bf 1}=\omega\ne 0$, a contradiction.
Thus $M$ is not a weak $V$-module.}
\end{rema}

\bigskip

\noindent {\small \sc Institute of Mathematics, Fudan University,
Shanghai 200433, China}

\noindent and

\noindent {\small \sc Department of Mathematics, Rutgers University,
Piscataway, NJ 08854}

\noindent (permanent address)

\noindent {\em E-mail address}: yzhuang@math.rutgers.edu

\vspace{1em}

\noindent {\small \sc Department of Mathematics, Rutgers University,
Piscataway, NJ 08854}

\noindent {\em E-mail address}: lepowsky@math.rutgers.edu

\vspace{1em}

\noindent {\small \sc Department of Mathematical Sciences, Rutgers
University, Camden, NJ 08102 and\\
Department of Mathematics, Harbin Normal University, Harbin, China}

\noindent {\em E-mail address}: hli@camden.rutgers.edu

\vspace{1em}

\noindent {\small \sc Department of Mathematics, Rutgers University,
Piscataway, NJ 08854}

\noindent {\em E-mail address}: linzhang@math.rutgers.edu


\begin{thebibliography}{ELM2}
\bibitem[B]{B}
R. E. Borcherds, Vertex algebras, Kac-Moody algebras, and the Monster,
{\em Proc. Natl. Acad. Sci. USA} {\bf 83} (1986), 3068--3071.

\bibitem[D]{D} C. Dong, Vertex algebras associated with even lattices,
{\em J.  Algebra} {\bf 160} (1993), 245--265.

\bibitem[DL]{DL} C. Dong and J. Lepowsky, {\em Generalized Vertex
Algebras and Relative Vertex Operators}, Progress in Math., Vol. 112,
Birkh\"{a}user, Boston, 1993.

\bibitem[DLM1]{DLM1} C. Dong, H.-S. Li and G. Mason, Vertex operator
algebras and associative algebras, {\em J. Algebra} {\bf 206} (1998),
67-96.

\bibitem[DLM2]{dlm-vpa}
C. Dong, H.-S. Li and G. Mason, Vertex Lie
algebra, vertex poisson algebras and vertex algebras, in: {\em Recent
Developments in Infinite-Dimensional Lie Algebras and Conformal Field
Theory}, Proceedings of an International Conference at University of
Virginia, May 2000, Contemporary Math. {\bf 297} (2002), 69-96.

\bibitem[DM]{DM} C. Dong and G. Mason, On quantum Galois theory, {\em
Duke Math. J.}  {\bf 86} (1997), 305--321.

\bibitem[FFR]{FFR}
A.~J. Feingold, I.~B. Frenkel and J.~F.~X. Ries, {\it Spinor
Construction of Vertex Operator Algebras, Triality, and
$E_{8}^{(1)}$}, Contemporary Math. {\bf 121} (1991).

\bibitem[FHL]{FHL}
I.~B. Frenkel, Y.-Z. Huang and J.~Lepowsky, On axiomatic approaches to
vertex operator algebras and modules, preprint, 1989; {\em Memoirs
Amer. Math. Soc.} {\bf 104}, 1993.

\bibitem[FLM]{FLM}
I.~B. Frenkel, J.~Lepowsky and A.~Meurman,
{\em Vertex Operator Algebras and the Monster},
Pure and Appl. Math., Vol. 134,  Academic Press,  Boston, 1988.

\bibitem[H1]{Hbook}
Y.-Z. Huang, {\em Two-dimensional Conformal Field Theory and Vertex
Operator Algebras}, Progress in Math., Vol. 148, Birkh\"{a}user, Boston, 1997.

\bibitem[H2]{tensor4}
Y.-Z. Huang, A theory of tensor products for module categories for a
vertex operator algebra, IV, {\em J. Pure Appl. Alg.} 100 (1995)
173--216.

\bibitem[HL1]{tensor1}
Y.-Z. Huang and J. Lepowsky, A theory of tensor products for module
categories for a vertex operator algebra, I, {\em Selecta Mathematica
(New Series)} {\bf 1} (1995), 699--756.

\bibitem[HL2]{tensor2}
Y.-Z. Huang and J. Lepowsky, A theory of tensor products for module
categories for a vertex operator algebra, II, {\em Selecta Mathematica
(New Series)} {\bf 1} (1995), 757--786.

\bibitem[HL3]{tensor3}
Y.-Z. Huang and J. Lepowsky, A theory of tensor
products for module categories for a vertex operator algebra, III,
{\em J. Pure Appl. Alg.} {\bf 100} (1995) 141--171.

\bibitem[HLZ1]{HLZ1} Y.-Z. Huang, J. Lepowsky and L. Zhang, A
logarithmic generalization of tensor product theory for modules for a
vertex operator algebra, arXiv:math.QA/0311235.

\bibitem[HLZ2]{HLZ2} Y.-Z. Huang, J. Lepowsky and L. Zhang,
Logarithmic tensor product theory for generalized modules for a
conformal vertex algebra, to appear.

\bibitem[LL]{LL} J. Lepowsky and H.-S. Li, {\em Introduction to Vertex
Operator Algebras and Their Representations}, Progress in Math.,
Vol. 227, Birkh\"auser, Boston, 2003.

\bibitem[Li1]{Li-form}
H.-S. Li, Symmetric invariant bilinear forms on vertex operator algebras,
{\em J. Pure Appl. Alg.} {\bf 96} (1994), 279-297.

\bibitem[Li2]{Li}
H.-S. Li, Representation theory and tensor product theory for vertex
operator algebras, Ph.D. thesis, Rutgers University, 1994.

\bibitem[Li3]{Li-local}
H.-S. Li, Local systems of vertex operators, vertex superalgebras and
modules, {\em J. Pure Appl. Alg.} {\bf 109} (1996), 143--195.

\bibitem[Li4]{Li-hom}
H.-S. Li, An analogue of the hom functor and a generalized nuclear
democracy theorem, {\em Duke Math. J.} {\bf 93} (1998), 73--114.

\bibitem[Mi]{Mi} A. Milas, Weak modules and logarithmic intertwining
operators for vertex operator algebras, in: {\em Recent Developments
in Infinite-Dimensional Lie Algebras and Conformal Field Theory},
Proceedings of an International Conference at University of Virginia,
May 2000, Contemporary Math. {\bf 297} (2002), 201--225.

\bibitem[MP]{MP} A. Meurman and M. Primc, Annihilating fields of
standard modules of $\widetilde{\mathfrak{s}\mathfrak{l}(2,\C)}$ and
combinatorial identities, Memoirs Amer. Math.  Soc. {\bf 652}, 1999.

\bibitem[P]{Primc} M. Primc, Vertex algebras generated by Lie
algebras, {\em J. Pure Applied Algebra} {\bf 135} (1999), 253--293.

\bibitem[Z]{Zhu}
Y. Zhu, Modular invariance of characters of vertex operator algebras,
{\em J. Amer. Math. Soc.} {\bf 9} (1996), 237--307.

\end{thebibliography}
\end{document}